%% file: main.tex
\documentclass[12pt,oneside,reqno]{amsart}
\input{preamble}

\begin{document} 
    
\input{title}


\input{abstract}
\maketitle

\bigskip
\bigskip

\setcounter{tocdepth}{2}
\tableofcontents

\input{body.tex}

\input{appendix.tex}
  \bigskip

\input{acknowledgment.tex}

\bibliographystyle{abbrvnat}
\bibliography{references}
\end{document}

%% file: preamble.tex
\usepackage{geometry}
\allowdisplaybreaks

\usepackage[utf8]{inputenc}
\usepackage[T1]{fontenc}
\usepackage[english]{babel}

\usepackage{lmodern}
\usepackage[nopatch=footnote]{microtype}

\usepackage{hyphenat}

\usepackage[dvipsnames,svgnames]{xcolor}

\definecolor{tocblue}{RGB}{0,158,196}

\usepackage{amsmath}
\usepackage{amsthm}
\usepackage{amssymb}
\usepackage{amsfonts}

\usepackage{mathtools}
\usepackage{nccmath}
\usepackage{empheq}

\usepackage{mathrsfs}
\usepackage{braket}
\usepackage{bm}

\usepackage{graphicx}

\usepackage{caption}
\usepackage{subcaption}

\usepackage{tikz}
\usepackage{pgfplots}

\pgfplotsset{compat=1.18}

\usepgfplotslibrary{fillbetween}

\usetikzlibrary{
    arrows,
    arrows.meta,
    positioning,
    patterns
}

\usepackage{enumitem}

\usepackage{verbatim}
\usepackage{comment}

\usepackage{fmtcount}
\usepackage{currvita}
\usepackage{csquotes}

\usepackage[numbers]{natbib}

\usepackage[
    colorlinks=true,
    linkcolor=Maroon,
    citecolor=Maroon,
    urlcolor=Maroon,
    linktoc=all,
    bookmarksopen=true,
    bookmarksnumbered=true
]{hyperref}

\usepackage{orcidlink}

\hypersetup{
    linkcolor=Maroon,
    citecolor=Maroon,
    urlcolor=Maroon,
    pdftitle={On the modified Zakharov-Kuznetsov equation on the cylinder: A trilinear Bourgain-space estimate and its application to local well-posedness},
    pdfauthor={Ali Mezher},
    pdfsubject={Modified Zakharov-Kuznetsov equation on the cylinder},
    pdfkeywords={modified Zakharov-Kuznetsov equation, Bourgain spaces, trilinear estimate, local well-posedness}
}

\newcommand{\assumref}[1]{%
    \hyperref[#1]{%
        \textcolor{Maroon}{Assumption~\ref*{#1}}%
    }%
}

\newcommand{\secref}[1]{%
    \hyperref[#1]{%
        \textcolor{Maroon}{Section~\ref*{#1}}%
    }%
}

\newcommand{\hypref}[1]{%
    \hyperref[#1]{%
        \textcolor{Maroon}{Hypothesis~\ref*{#1}}%
    }%
}

\newcommand{\appref}[1]{%
    \hyperref[#1]{%
        \textcolor{Maroon}{Appendix~\ref*{#1}}%
    }%
}

\newcommand{\apprefs}[2]{%
    \hyperref[#1]{%
        \textcolor{Maroon}{Appendices~\ref*{#1}--\ref*{#2}}%
    }%
}

\newcommand{\appsubref}[1]{%
    \hyperref[#1]{%
        \textcolor{Maroon}{Appendix~\ref*{#1}}%
    }%
}

\newcommand{\thmref}[1]{%
    \hyperref[#1]{%
        \textcolor{Maroon}{Theorem~\ref*{#1}}%
    }%
}

\newcommand{\lemref}[1]{%
    \hyperref[#1]{%
        \textcolor{Maroon}{Lemma~\ref*{#1}}%
    }%
}

\newcommand{\propref}[1]{%
    \hyperref[#1]{%
        \textcolor{Maroon}{Proposition~\ref*{#1}}%
    }%
}

\newcommand{\corref}[1]{%
    \hyperref[#1]{%
        \textcolor{Maroon}{Corollary~\ref*{#1}}%
    }%
}

\newcommand{\remref}[1]{%
    \hyperref[#1]{%
        \textcolor{Maroon}{Remark~\ref*{#1}}%
    }%
}

\newcommand{\figref}[1]{%
    \hyperref[#1]{%
        \textcolor{Maroon}{Figure~\ref*{#1}}%
    }%
}

\newcommand{\eqrefmaroon}[1]{%
    \hyperref[#1]{%
        \textcolor{Maroon}{(\ref*{#1})}%
    }%
}

\numberwithin{equation}{section}

\theoremstyle{plain}

\newtheorem{theorem}{Theorem}[section]

\newtheorem{definition}[theorem]{Definition}

\newtheorem{lemma}[theorem]{Lemma}

\newtheorem{remark}[theorem]{Remark}

\newtheorem*{remark*}{Remark}

\makeatletter

\def\namedlabel#1#2{%
    \begingroup
    #2%
    \def\@currentlabel{#2}%
    \phantomsection
    \label{#1}%
    \endgroup
}

\makeatother

\makeatletter

\def\paragraph{%
    \@startsection{paragraph}{4}%
    \z@
    \z@
    {-\fontdimen2\font}%
    {\normalfont\bfseries}%
}

\makeatother

\makeatletter

\def\@tocline#1#2#3#4#5#6#7{%
    \relax
    \ifnum #1>\c@tocdepth
    \else
        \par
        \addpenalty\@secpenalty
        \addvspace{#2}%
        \begingroup
        \hyphenpenalty\@M
        \@ifempty{#4}{%
            \@tempdima
            \csname r@tocindent\number#1\endcsname\relax
        }{%
            \@tempdima#4\relax
        }%
        \parindent\z@
        \leftskip#3\relax
        \advance\leftskip\@tempdima\relax
        \rightskip\@pnumwidth plus4em
        \parfillskip-\@pnumwidth
        #5\leavevmode
        \hskip-\@tempdima
        \ifcase#1
        \or
        \or
            \hskip2em
        \or
            \hskip3em
        \else
            \hskip4em
        \fi
        #6\nobreak\relax
        \hfill
        \hbox to\@pnumwidth{%
            \@tocpagenum{#7}%
        }%
        \par
        \nobreak
        \endgroup
    \fi
}

\makeatother

\def\C{{\mathbb C}}
\def\N{{\mathbb N}}

\def\R{{\mathbb R}}
\def\Z{{\mathbb Z}}

\def\T{{\mathbb T}}

\newcommand{\supp}{\operatorname{supp}}

\def\cF{{\mathcal F}}

\DeclareRobustCommand{\rchi}{%
    \raisebox{0.3ex}{$\chi$}%
}

\def\1{{\mathbf 1}}

% ...
% ...

\newcommand{\abs}[1]{%
    \left|#1\right|%
}

\newcommand{\norm}[1]{%
    \left\lVert #1\right\rVert%
}

\newcommand{\inp}[1]{%
    \left\langle #1\right\rangle
}

\newcommand{\prc}[1]{%
    \left(#1\right)
}

\newcommand{\brs}[1]{%
    \left\{#1\right\}
}

\def\pt{\partial}

\newcommand{\mcapinn}[2]{%
    \vcenter{%
        \hbox{%
            \scalebox{0.7}{%
                $\mathsurround=0pt
                \ifx\displaystyle#1\textstyle
                \else
                    #1
                \fi
                \bigcap$
            }%
        }%
    }%
}

\newcommand{\mcupinn}[2]{%
    \vcenter{%
        \hbox{%
            \scalebox{0.7}{%
                $\mathsurround=0pt
                \ifx\displaystyle#1\textstyle
                \else
                    #1
                \fi
                \bigcup$
            }%
        }%
    }%
}

\newcommand{\RomanNumeralCaps}[1]{%
    \uppercase\expandafter{%
        \romannumeral #1\relax
    }%
}

\makeatletter

\newcommand{\rom}[1]{%
    \expandafter\@slowromancap
    \romannumeral #1@%
}

\makeatother

\newlength\mylen
\newlist{mycases}{enumerate}{1}

\setlist[mycases,1]{
    label=\textbf{Case~\arabic*:},
    labelwidth=*,
    leftmargin=*,
    align=right
}

\newlength\mylenn
\newlist{mycasesproof}{enumerate}{1}

\setlist[mycasesproof,1]{
    label=\textbf{Proof of Case~\arabic*:},
    labelwidth=*,
    leftmargin=*,
    align=right
}

\newcommand{\ds}{\displaystyle}

\newcommand{\Lpn}[2]{%
    \left\lVert #1\right\rVert_{L^{#2}}%
}

\newcommand{\Lptxy}[3]{%
    \left\lVert #1\right\rVert_{L^{#2}_{#3}}%
}

\newcommand{\Hn}[2]{%
    \left\lVert #1\right\rVert_{H^{#2}}%
}

\newcommand{\Hnxy}[3]{%
    \left\lVert #1\right\rVert_{H^{#2}_{#3}}%
}

\newcommand{\xnorm}[3]{%
    \left\lVert #1\right\rVert_{X^{{#2},{#3}}}%
}

\newcommand{\xtnorm}[3]{%
    \left\lVert #1\right\rVert_{X^{{#2},{#3}}_T}%
}

\newcommand{\rt}{\R\times\T}

\newcommand{\rrt}{\R\times\R\times\T}
\newcommand{\rrz}{\R\times\R\times\Z}

\newcommand{\F}{\mathcal{F}}
\newcommand{\Fh}[1]{\widehat{#1}}

\newcommand{\Ft}{\mathcal{F}_t}

\newcommand{\Fxy}{\mathcal{F}_{xy}}
\newcommand{\IF}{{\mathcal{F}^{-1}}}

\newcommand{\IFxy}{\mathcal{F}_{xy}^{-1}}

\newcommand{\PQ}[2]{%
    P_{N_{#2}}Q_{L_{#2}}#1_{#2}%
}

\newcommand{\FPQ}[2]{%
    \widehat{P_{N_{#2}}Q_{L_{#2}}#1_{#2}}%
}

\newcommand{\FPQwV}[2]{%
    \widehat{%
        P_{N_{#2}}Q_{L_{#2}}#1_{#2}%
    }%
    (\tau_{#2},\xi_{#2},q_{#2})%
}

\newcommand{\Ifr}[4]{%
    I_{{#1}{#2}{#3}\to{#4}}%
}

\newcommand{\Ns}{{N,N_1,N_2,N_3}}
\newcommand{\Ls}{{L,L_1,L_2,L_3}}

\newcommand{\prdd}{%
    \prod_{j=1}^{3}%
}

\newcommand{\sg}[1]{%
    e^{#1\partial_x\Delta}%
}

%% file: title.tex
\title[On the modified Zakharov--Kuznetsov equation on the cylinder]
{On the modified Zakharov--Kuznetsov equation on the cylinder:
A trilinear Bourgain-space estimate and its application to local well-posedness}
    
\author{Ali Mezher\orcidlink{0000-0002-9215-1804}}

\address{Department of Mathematics, University of Texas at Austin, Austin, TX 78712, USA}
\email{alimezher@utexas.edu}

%% file: abstract.tex
\begin{abstract}
We study the modified Zakharov--Kuznetsov equation on $\R\times\T$. We establish a trilinear Bourgain-space estimate on the full time line for the derivative cubic term, with input exponent $\frac12+\delta$ and output exponent $-\frac12+3\delta$, and use it to recover local well-posedness in $H^s(\R\times\T)$ for $s>1$. After time localization, this output exponent yields a factor $T^{2\delta}$ in the contraction argument. Lower Sobolev regularity thresholds are already known; the purpose here is to give a self-contained dyadic proof of this different modulation balance.
\end{abstract}

%% file: body.tex
\section{Introduction}

The modified Zakharov--Kuznetsov equation is a multidimensional dispersive
model of Korteweg--de Vries type in which the quadratic nonlinearity of the
Zakharov--Kuznetsov equation is replaced by a cubic one. In this work, we
consider the initial-value problem on the semi-periodic domain
$\R\times\T$,
\begin{equation}
\label{eq:mZK-cylinder-intro}
\begin{cases}
    \partial_t u+\partial_x\Delta u+\partial_x(u^3)=0,\vspace{ 0.6 em}
    \\
    u(0,x,y)=u_0(x,y),
\end{cases}
\qquad
(x,y)\in\R\times\T.
\end{equation}
The semi-periodic geometry combines a continuous longitudinal frequency with
a discrete transverse frequency and consequently produces a resonance
structure that differs substantially from that of the purely Euclidean
problem. At the Fourier level, the linear evolution associated with
\eqref{eq:mZK-cylinder-intro} is governed by the dispersion relation
\begin{equation}
\label{eq:mZK-dispersion-intro}
    \omega(\xi,q)=\xi^3+\xi q^2,
    \qquad
    (\xi,q)\in\R\times\Z.
\end{equation}
Thus, the analysis of the derivative cubic nonlinearity requires simultaneous
control of spatial frequencies, modulation variables, and Sobolev weights,
while retaining the interaction between the continuous frequency $\xi$ and
the discrete transverse mode $q$. The purpose of the present paper is to
establish and apply a Bourgain-space estimate on the full timeline, with input
modulation exponent $\frac12+\delta$ and output modulation exponent
$-\frac12+3\delta$. Direct trilinear Bourgain-space estimates for the same
cubic derivative term already appears in
\cite{farah2024note,nowickikoth2026strichartz}; the precise relation between
those estimates and the estimate proved here are described below.

The Zakharov--Kuznetsov equation was introduced by Zakharov and Kuznetsov as
a multidimensional model for the propagation of nonlinear ion-acoustic waves
in magnetized plasmas \cite{zakharov1974three}; related derivations and
physical aspects can be found in
\cite{laedke1982nonlinear,infeld2000nonlinear}. A rigorous derivation from the
Euler--Poisson system in the presence of a magnetic field was subsequently
obtained by Lannes, Linares, and Saut \cite{lannes2013cauchy}. Modified
versions of the equation arise when the leading nonlinear interaction is
cubic rather than quadratic. In particular, Munro and Parkes derived a
modified Zakharov--Kuznetsov equation and studied the stability of its
solutions \cite{munro1999derivation}; further dynamical aspects of the
modified equation were considered in \cite{sipcic2000lump}. From the
analytical point of view, this change in the nonlinear structure leads
naturally from the bilinear estimates associated with the quadratic ZK
nonlinearity to multilinear estimates capable of controlling a derivative
cubic interaction.

For sufficiently smooth solutions of \eqref{eq:mZK-cylinder-intro}, the mass
and Hamiltonian energy are conserved
\begin{equation}
\begin{aligned}
    M(u(t))
    &:=
    \int_{\R\times\T}
    |u(t,x,y)|^2\,dx\,dy,
    \\*
    E(u(t))
    &:=
    \frac12
    \int_{\R\times\T}
    |\nabla u(t,x,y)|^2\,dx\,dy
    -
    \frac14
    \int_{\R\times\T}
    |u(t,x,y)|^4\,dx\,dy.
\end{aligned}
\end{equation}
Consequently, the energy space $H^1(\R\times\T)$ plays a distinguished role in the dynamics. On $\R^d$, the modified equation possesses the scaling
\begin{equation*}
    u_\lambda(t,x)
    =
    \lambda u(\lambda^3t,\lambda x),
\end{equation*}
under which
\begin{equation*}
    \norm{u_\lambda(0)}_{\dot H^s(\R^d)}
    =
    \lambda^{s+1-\frac d2}
    \norm{u_0}_{\dot H^s(\R^d)}.
\end{equation*}
Hence the Euclidean scaling-critical Sobolev index is
\begin{equation*}
    s_c=\frac d2-1.
\end{equation*}
The isotropic Euclidean scaling, however, does not preserve
$\R\times\T$, since it changes the period in the transverse variable. Therefore, on the cylinder, scaling alone does not determine the expected
regularity threshold. The mixed continuous-discrete resonance geometry
becomes an essential part of the low-regularity analysis.

The well-posedness theory for the modified equation on Euclidean spaces has
developed through a sequence of increasingly refined dispersive and
multilinear estimates. Biagioni and Linares established local well-posedness
for the two-dimensional mZK equation in $H^1(\R^2)$
\cite{biagioni2003well}. Linares and Pastor subsequently proved local
well-posedness in $H^s(\R^2)$ for $s>\frac34$ and showed that well-posedness
cannot hold in the scaling-critical space $L^2(\R^2)$
\cite{linares2009well}. They also obtained global results for the
two-dimensional generalized ZK equation; in the cubic case, under the
appropriate mass restriction in the focusing regime, they proved global
well-posedness for $s>\frac{53}{63}$
\cite{linares2011local}. Ribaud and Vento lowered the local threshold for
the two-dimensional cubic generalized ZK equation to
$s>\frac14$ \cite{ribaud2012note}, while related local, global, and
scattering results for generalized ZK equations were developed in
\cite{farah2012note}. Kinoshita subsequently reached the endpoint
$s=\frac14$ for the two-dimensional modified equation and established the
sharpness of this regularity in terms of the smoothness of the flow map
\cite{kinoshita2022well}. In three dimensions, Gr\"unrock proved local
well-posedness for $s>\frac12$ \cite{grunrock2014remark}. Global
low-regularity theory in two dimensions was
further developed by Bhattacharya, Farah, and Roudenko using the $I$-method
\cite{bhattacharya2020global}. More recently, Correia and Kinoshita introduced
anisotropic Sobolev-type spaces adapted to the dispersive geometry of the
two-dimensional modified equation and established sharp local and small-data
global scattering results in that framework
\cite{correia2026global}.

The semi-periodic problem has a different history. For the original
Zakharov--Kuznetsov equation on $\R\times\T$, Linares, Pastor, and Saut
initiated the well-posedness theory on the cylinder and obtained local
well-posedness in $H^s(\R\times\T)$ for $s>\frac32$
\cite{linares2010well}. Molinet and Pilod subsequently developed bilinear
Strichartz estimates adapted to the ZK dispersion relation and obtained,
among other results, global well-posedness in the energy space
$H^1(\R\times\T)$ \cite{molinet2015bilinear}. Osawa extended the local
theory below the energy level, proving local well-posedness for
$\frac{9}{10}<s<1$ \cite{osawa2022local}. More recently, Cao-Labora refined
the analysis of the resonant set and proved local well-posedness for arbitrary
$H^s(\R\times\T)$ data when $s>\frac34$ and, under an additional
low-frequency condition on the initial data, when $s>\frac12$
\cite{cao2025low}. These developments illustrate a characteristic feature of
semi-periodic dispersive equations: the mixture of a continuous and a
discrete Fourier variable creates resonant configurations that require
estimates specifically adapted to the geometry of the cylinder.

For higher-order nonlinearities on $\R\times\T$, Farah and Molinet developed
a local and global theory for the generalized Zakharov--Kuznetsov equation
\cite{farah2024note}. In the modified case, corresponding to the cubic
nonlinearity, their Bourgain-space argument yields local well-posedness for
$s>\frac56$, together with global results in the energy space under
appropriate smallness assumptions on the initial data. More recently,
Nowicki-Koth established an almost optimal linear $L^4$ Strichartz estimate
and a family of bilinear refinements adapted to the mixed
continuous-discrete Fourier geometry. As an application, the local
well-posedness threshold for the semi-periodic modified
Zakharov--Kuznetsov equation was lowered to $ s>{11}/{24}$; see \cite{nowickikoth2026strichartz}.

In subsequent work, these estimates were combined with the $I$-method to obtain global well-posedness for initial data satisfying
\begin{equation}
    \norm{u_0}_{L^2(\R\times\T)}<\norm{Q_2}_{L^2(\R^2)},
\end{equation}
where $Q_2$ denotes
the two-dimensional ground state, whenever $s>{36}/{49}$; see \cite{nowickikoth2026global}. 

Thus, the currently available
low-regularity theory on the cylinder exploits a delicate combination of
linear Strichartz estimates, bilinear smoothing, resonance analysis, and the
arithmetic structure of the transverse frequencies.

Direct trilinear Bourgain-space estimates for the cubic derivative term are
already available in the literature. In the modified case $k=2$,
\cite[Proposition~3.3]{farah2024note} implies, after placing all three inputs
at the same spatial regularity, that for every $s>\frac56$ and sufficiently
small $\varepsilon>0$,
\begin{equation}
\label{eq:farah-molinet-trilinear-estimate}
    \left\|
        \partial_x(u_1u_2u_3)
    \right\|_{X^{s,-\frac12+2\varepsilon}}
    \lesssim
    \prod_{j=1}^{3}
    \norm{u_j}_{X^{s,\frac12+\varepsilon}}.
\end{equation}
More recently, \cite[Proposition~4.3]{nowickikoth2026strichartz} established
the same symmetric modulation balance for $s>\frac{11}{24}$ under a
time-localization assumption on the inputs.

The estimate proved in the present paper has a different modulation balance
and is stated on the full timeline. More precisely, for every $s>1$, there
exists $\delta=\delta(s)>0$ such that
\begin{equation}
\label{eq:intro-trilinear-estimate}
    \left\|
        \partial_x(u_1u_2u_3)
    \right\|_{X^{s,-\frac12+3\delta}}
    \lesssim
    \prod_{j=1}^{3}
    \norm{u_j}_{X^{s,\frac12+\delta}}.
\end{equation}
For the same positive parameter, the target space in
\eqref{eq:intro-trilinear-estimate} has a higher modulation exponent than the
target space in \eqref{eq:farah-molinet-trilinear-estimate}. This exact balance
does not follow from the preceding estimates merely by monotonicity of
Bourgain norms: applying an estimate with parameter $\varepsilon$ to inputs
known only to belong to $X^{s,\frac12+\delta}$ requires
$\varepsilon\leq\delta$, whereas control of the target norm in
\eqref{eq:intro-trilinear-estimate} would require
$2\varepsilon\geq3\delta$. On the other hand, the condition $s>1$ is more
restrictive than the spatial regularity ranges in the cited works.

The proof of \eqref{eq:intro-trilinear-estimate} uses a dyadic decomposition
of both spatial frequencies and modulations. After ordering the input
frequencies, we consider the six regimes low--low--low $\to$ low,
high--low--low $\to$ high, high--high--low $\to$ low,
high--high--low $\to$ high, high--high--high $\to$ low, and
high--high--high $\to$ high. Frequency-separated interactions are controlled
by first-derivative bounds for the resonance function. In the
high--high--high $\to$ high regime, we retain a dyadic pair-frequency
parameter $K$ and use a quadratic-phase sublevel-set estimate. The underlying
$K^{-\frac14}$ bilinear gain is the pair-frequency-localized form of the
high-output-frequency estimate appearing in the proof of
\cite[Proposition~3.1]{farah2024note}. In the present argument, that
$K$-dependence is retained, and the bilinear estimate is applied to two paired
convolutions in the dual formulation of the trilinear bound.

The condition $s>1$ is the threshold at which the dyadic summations in this
particular proof close; no claim of optimality is made. Combined with the
inhomogeneous Bourgain-space estimate, the output exponent
$-\frac12+3\delta$ produces the factor $T^{2\delta}$ used in the contraction
argument below.

As an application of \eqref{eq:intro-trilinear-estimate}, together with the
standard homogeneous and inhomogeneous linear estimates in Bourgain spaces,
we obtain local well-posedness of \eqref{eq:mZK-cylinder-intro} in
$H^s(\R\times\T)$ for every $s>1$. The solution is constructed in the
restricted Bourgain space $X_T^{s,\frac12+\delta}$, which embeds continuously
into $C([0,T];H^s(\R\times\T))$. The same framework yields uniqueness in the
full restricted Bourgain class. Finally, the analytic structure of the
nonlinear mapping, together with the implicit function theorem in Banach
spaces, gives real-analytic dependence of the solution on the initial data.
We now state the main result.

\begin{theorem}
\label{MainResult}
Assume that $s>1$. Then there exists $\delta=\delta(s)>0$ such that, for
every $u_0\in H^s(\R\times\T)$, there exists
\begin{equation*}
    0<T=T\prc{\norm{u_0}_{H^s(\R\times\T)}}
\end{equation*}
for which the initial-value problem
\eqref{eq:mZK-cylinder-intro} admits a unique solution
\begin{equation*}
    u\in
    C\prc{[0,T];H^s(\R\times\T)}
    \cap
    X_T^{s,\frac12+\delta}.
\end{equation*}
Moreover, for every $T'\in(0,T)$, there exists $r>0$ such that, with
\begin{equation*}
    B_r(u_0)
    :=
    \brs{
        v_0\in H^s(\R\times\T):
        \norm{v_0-u_0}_{H^s(\R\times\T)}<r
    },
\end{equation*}
the data-to-solution map
\begin{equation*}
\begin{aligned}
    \mathcal S_r:
    B_r(u_0)
    &\longrightarrow
    C\prc{[0,T'];H^s(\R\times\T)}
    \cap
    X_{T'}^{s,\frac12+\delta},
    \\
    v_0
    &\longmapsto
    v
\end{aligned}
\end{equation*}
is real-analytic.
\end{theorem}

The restricted Bourgain space $X_T^{s,b}$ is defined in \secref{sec:function-space}.

\section{Notations and Function Spaces}
\label{sec:function-space}
\subsection{Notations}
\label{sec:notations}
For any positive numbers $a$ and $b$, the notation $a \lesssim b$ means that there exists a positive constant $c$ such that $a \leq cb$. Moreover, we write $a \lesssim_d b$ when $a\leq c b$ and $c$ is a function of $d$, namely $c=c(d)$. We also write $a \sim b$ when $a \lesssim b$ and $b \lesssim a$. If $A$ and $B$ are two positive numbers, we use the notation $A \land B = \min(A, B)$ and $A \lor B = \max(A, B)$. Finally, $|\cdot|$ denotes the Lebesgue measure for measurable subsets of $\R^m$, the counting measure for subsets of $\Z^n$, and the corresponding product measure on $\R^m\times\Z^n$, where $m,n\in\Z^+$. For $(\xi,q)\in\R\times\Z$, we use the notation
\begin{equation}
    \abs{(\xi,q)}
    :=
    \sqrt{3\xi^2+q^2}.
\end{equation}
Since
\begin{equation}
    \omega(\xi,q)
    :=
    \xi^3+\xi q^2,
\end{equation}
we have
\begin{equation}
    \pt_\xi\omega(\xi,q)
    =
    3\xi^2+q^2
    =
    \abs{(\xi,q)}^2.
\end{equation}
The function $\omega$ is the dispersion relation associated with the
linear part of the ZK and mZK equations.

We write $\cF(u)$, or $\hat u$, for the space--time Fourier transform of $u$,
and $\Fxy(u)$ and $\Ft(u)$ for its spatial and temporal Fourier transforms,
respectively. For $s\in\R$, we define the Bessel and Riesz potentials of
order $s$, $J^s$ and $D^s$, by

\begin{align}
    J^s u &=\IFxy((1+ |(\xi,q)|^2)^\frac{s}{2} \Fxy(u) ),\\
    D^s u &=\IFxy( |(\xi,q)|^s \Fxy(u) ).
\end{align}
We fix an even function $\rchi\in C_0^\infty(\R)$ such that
$0\leq\rchi\leq1$, the restriction of $\rchi$ to $[0,\infty)$ is
nonincreasing,
\begin{equation}
    \rchi|_{[-5/4,5/4]}=1,
    \qquad
    \supp(\rchi)\subset
    \left[-\frac85,\frac85\right].
\end{equation}
For $k\in\Z^+$, define
\begin{equation*}
    \phi(r):=\rchi(r)-\rchi(2r),
    \qquad
    \phi_{2^k}(\xi,q)
    :=
    \phi\big(2^{-k}\abs{(\xi,q)}\big),
\end{equation*}
and
\begin{equation*}
    \psi_{2^k}(\tau,\xi,q)
    :=
    \phi\big(2^{-k}(\tau-\omega(\xi,q))\big).
\end{equation*}
Moreover, let
\begin{equation*}
    \phi_1(\xi,q):=\rchi\big(\abs{(\xi,q)}\big),
    \qquad
    \psi_1(\tau,\xi,q):=\rchi(\tau-\omega(\xi,q)).
\end{equation*}
Summations over $N$, $L$, or $K$ are understood to be dyadic, with
$N,L,K\geq1$. For a dyadic number $M\geq1$, set
\begin{equation*}
    I_1:=\left[0,\frac85\right],
    \qquad
    I_M:=\left[\frac58M,\frac85M\right]
    \quad\text{for }M\geq2.
\end{equation*}
Then
\begin{equation*}
    \sum_N\phi_N(\xi,q)=1,
    \qquad
    \supp(\phi_N)
    \subset
    \left\{(\xi,q):\abs{(\xi,q)}\in I_N\right\},
\end{equation*}
and
\begin{equation*}
    \sum_L\psi_L(\tau,\xi,q)=1,
    \qquad
    \supp(\psi_L)
    \subset
    \left\{(\tau,\xi,q):
    \abs{\tau-\omega(\xi,q)}\in I_L\right\}.
\end{equation*}
Define the spatial Littlewood--Paley multiplier $P_N$ and the modulation
multiplier $Q_L$ by

$$P_N (u) = \IFxy{ \prc{\phi_N \Fxy(u)}}, \hspace{1cm} Q_L (u) = \IF{ \prc{\psi_L \F(u)}}.$$
Finally, we denote by $\brs{e^{- t \pt_x \Delta}}_{t\in \R}$ the free group associated with the linear part of the equation \eqref{eq:mZK-cylinder-intro}; that is,

$$\Fxy{\prc{e^{-t \pt_x \Delta} \varphi}}= e^{i t \omega(\xi,q)} \Fxy(\varphi),$$
where $\varphi \in L^2 (\R \times \T)$.

\smallskip
The resonance function is denoted by $\mathcal{H}$ and defined by

\begin{equation}
    \begin{aligned}
\label{ResonanceFunction}
    \mathcal{H}(\xi_1,q_1,\xi_2,q_2,\xi_3, q_3)=& \omega(\xi_1 + \xi_2 +\xi_3,q_1 +q_2+q_3)\\
    & - \omega(\xi_1,q_1)-\omega( \xi_2,q_2) -\omega(\xi_3,q_3).
\end{aligned}
\end{equation}

\subsection{Function spaces}
For $1\leq p\leq\infty$, $L^p(\R\times\T)$ is the usual Lebesgue space
with norm $\Lpn{\cdot}{p}$. For $s\in\R$, the real-valued Sobolev space
$H^s(\R\times\T)$ has norm $\Hn{u}{s}=\Lpn{J^su}{2}$. For $T>0$,
$1\leq p<\infty$, and a Banach space $B$ of functions on $\R\times\T$,
define the mixed space-time norms by
\begin{equation}
\label{eq:mixed-space-time-norms}
\begin{aligned}
\norm{u}_{L_T^pB_{xy}}
&:=
\left(
    \int_0^T
    \norm{u(t,\cdot,\cdot)}_B^p
    \,dt
\right)^{\frac1p},
\\
\norm{u}_{L_t^pB_{xy}}
&:=
\left(
    \int_{\R}
    \norm{u(t,\cdot,\cdot)}_B^p
    \,dt
\right)^{\frac1p}.
\end{aligned}
\end{equation}
For $p=\infty$, we use the corresponding essential-supremum norms. For
$s,b\in\R$, we denote by $X^{s,b}$
the Bourgain space associated with the linear part of
\eqref{eq:mZK-cylinder-intro}; see \cite{bourgain1993fourier}.

$$X^{s,b}=\brs{ u=u(t,x,y)\in\mathcal{S}'(\rrt); \xnorm{u}{s}{b} <\infty},$$
where
$$\xnorm{u}{s}{b} = \norm{\inp{\tau -\omega(\xi,q)}^{b} \inp{\abs{(\xi,q)}}^{s} \F(u)(\tau,\xi,q)}_{  L^2_{\tau \xi} \ell_q^2},$$
and $\inp{x}:=(1+|x|^2)^{1/2}$. Moreover, for $T>0$, define
\begin{equation*}
\begin{aligned}
X_T^{s,b}:=
\Big\{u=u(t,x,y)\in\mathcal{S}'([0,T]\times\rt):{}
&\exists v\in X^{s,b}\text{ such that}\\
&v|_{(0,T)\times\R\times\T}=u
\Big\}.
\end{aligned}
\end{equation*}
We equip $X_T^{s,b}$ with the norm
$$\norm{u}_{X_T^{s,b}} = \inf \brs{\norm{v}_{X^{s,b}}: v\in X^{s,b}, v|_{(0,T) \times \R \times \T}=u}.$$

\section{The trilinear estimate on the cylinder}
We use the following elementary one-dimensional sublevel-set estimate; see
also \cite[Lemma~3.8]{molinet2015bilinear}.

\begin{lemma}[One-dimensional sublevel-set estimate]
\label{Lem:Lemma3.8}
Let $I,J\subset\R$ be intervals and let $f\in C^1(J;\R)$. If
\begin{equation*}
    \lambda:=\inf_{x\in J}\abs{f'(x)}>0,
\end{equation*}
then
\begin{equation}
\label{eq:one-dimensional-sublevel-estimate}
    \abs{\brs{x\in J:f(x)\in I}}
    \leq
    \frac{\abs{I}}{\lambda}.
\end{equation}
\end{lemma}

\begin{proof}
Since $f'$ is continuous and does not vanish on the interval $J$, it has a fixed sign there. Thus $f$ is one-to-one on $J$, and the estimate follows
from the change-of-variables formula.
\end{proof}

\begin{lemma}
\label{Lem:TheTrilinearEstimateIn-L2}
Let $N_i,L_i$, $i=1,2,3$, be dyadic numbers in
$\brs{2^k:k\in\Z^+\cup\brs{0}}$, and assume that
$u_i(t,x,y)\in L^2(\rrt)$. We denote by $N_{\min}$ and $N_{\mathrm{med}}$ the
minimum and median, respectively, of $N_1,N_2,N_3$, and define
$L_{\min}$, $L_{\mathrm{med}}$, and $L_{\max}$ analogously. Then
\begin{equation}
\label{TheTrivialTrilinearEstimateIn-L2}
    \Lpn{\prod_{j=1}^3 {\PQ{u}{j}}}{2} \lesssim L_{\min}^\frac{1}{2} L_{\mathrm{med}}^\frac{1}{2} N_{\min} N_{\mathrm{med}} \prod_{j=1}^3 \Lpn{\PQ{u}{j}}{2}.
\end{equation}

Moreover,

\begin{enumerate}
\item[(\rom{1})] If $N_1 \geq 4 N_3$ or $N_3 \geq 4 N_1$, then
    \begin{equation}
    \label{eq:trilinear-L2-N1-N3}
\Lpn{\prod_{j=1}^3 {\PQ{u}{j}}}{2} \lesssim (N_{\min}N_{\mathrm{med}}L_1L_2L_3)^\frac{1}{2} \frac{N_2^\frac{1}{2}}{N_1 \lor N_3} \prod_{j=1}^3 \Lpn{\PQ{u}{j}}{2},
\end{equation}

\item[(\rom{2})] If $N_2 \geq 4 N_3$ or $N_3 \geq 4 N_2$, then 
 \begin{equation}
    \label{eq:trilinear-L2-N2-N3}
\Lpn{\prod_{j=1}^3 {\PQ{u}{j}}}{2} \lesssim (N_{\min}N_{\mathrm{med}}L_1L_2L_3)^\frac{1}{2} \frac{N_1^\frac{1}{2}}{N_2 \lor N_3} \prod_{j=1}^3 \Lpn{\PQ{u}{j}}{2},
\end{equation}

\item[(\rom{3})] If $N_1 \geq 4 N_2$ or $N_2 \geq 4 N_1$, then 
 \begin{equation}
    \label{eq:trilinear-L2-N1-N2}
\Lpn{\prod_{j=1}^3 {\PQ{u}{j}}}{2} \lesssim (N_{\min}N_{\mathrm{med}}L_1L_2L_3)^\frac{1}{2} \frac{N_3^\frac{1}{2}}{N_1 \lor N_2} \prod_{j=1}^3 \Lpn{\PQ{u}{j}}{2},
\end{equation}
    
\end{enumerate}

\end{lemma}

\begin{proof}
    By Plancherel's identity, we have
\begin{align*}
\Lptxy{
    \prod_{j=1}^3\PQ{u}{j}
}{2}{txy}
&=
\Lptxy{
    \mathcal F
    \left(
        \prod_{j=1}^3\PQ{u}{j}
    \right)
}{2}{\tau\xi q}
\\
&=
\left\|
    \int_{\R^4}
    \sum_{(q_1,q_2)\in\Z^2}
    \prod_{j=1}^3
    \FPQwV{u}{j}
    \,d\nu
\right\|_{L^2_{\tau,\xi,q}},
\end{align*}
where $\circledast_{j=1}^n u_j:=u_1\ast\cdots\ast u_n$, $\ast$ denotes
convolution, $d\nu:=d\tau_1\,d\tau_2\,d\xi_1\,d\xi_2$,
$\tau=\tau_1+\tau_2+\tau_3$, $\xi=\xi_1+\xi_2+\xi_3$, and
$q=q_1+q_2+q_3$. Applying Cauchy--Schwarz to the integral inside the norm and
using Fubini's theorem, the preceding expression is bounded by

\begin{align*}  
 &\Lptxy{\prc{\int_{\R^4}  \ds\sum_{(q_1,q_2) \in \Z^2} \mathbf 1_{A_{(\tau, \xi, q)} } d\nu}^\frac{1}{2} \prc{\int_{\R^4} \ds\sum_{(q_1,q_2) \in \Z^2} \abs{\prdd \FPQwV{u}{j}}^2 d\nu}^\frac{1}{2} }{2}{\tau \xi q}\\
& \leq \sup_{\tau,\xi,q}\prc{ \int_{\R^4} \sum_{(q_1,q_2) \in \Z^2} {\mathbf 1_{A_{(\tau, \xi, q)} }}  d\nu}^\frac{1}{2} \\
&\qquad \qquad \qquad  \times  \Lptxy{ \prc{\int_{\R^4} \ds\sum_{(q_1,q_2) \in \Z^2} \abs{\prdd \FPQwV{u}{j}}^2 d \nu}^\frac{1}{2} }{2}{\tau \xi q}\\
& \leq  \sup_{(\tau,\xi,q)}\abs{A_{(\tau,\xi,q)}}^\frac{1}{2}\\
&\times \Bigg(\int_{\R^2}\sum_{q} \abs{\prc{\int_{\R^4} \sum_{(q_1,q_2) \in \Z^2} \abs{\prdd \FPQwV{u}{j}}^2 d \nu}}d\tau d\xi\Bigg)^\frac{1}{2} \\
& \leq \sup_{(\tau,\xi,q)}\abs{A_{(\tau,\xi,q)}}^\frac{1}{2} \Bigg(\int_{\R^4}\sum_{(q_1,q_2) \in \Z^2}   \abs{\prod_{j=1}^2\FPQwV{u}{j}}^2\\
&\qquad \qquad \qquad \qquad \qquad \qquad \qquad  \times \prc{\int_{\R^2} \sum_q \abs{\FPQ{u}{3}}^2 d\tau d\xi} d\nu\Bigg)^\frac{1}{2}\\
&= \sup_{(\tau,\xi,q)}\abs{A_{(\tau,\xi,q)}}^\frac{1}{2}   \Lptxy{\FPQ{u}{3}}{2}{\tau \xi q}\\
&\qquad \qquad \qquad \qquad  \times \prc{\int_{\R^4}\sum_{(q_1,q_2) \in \Z^2}   \abs{\prod_{j=1}^2\FPQwV{u}{j}}^2d\nu}^\frac{1}{2}\\
&= \sup_{(\tau,\xi,q)}\abs{A_{(\tau,\xi,q)}}^\frac{1}{2} \prod_{j=1}^3 \Lptxy{\PQ{u}{j}}{2}{t x y},
\end{align*}
where 
\begin{equation}
\begin{aligned}
\label{TheSetA}
A_{(\tau,\xi,q)} = &\{ ((\tau_1,\xi_1,q_1),(\tau_2,\xi_2,q_2)) \in \R^4 \times \Z^2: \abs{(\xi_1,q_1)} \in I_{N_1},\\
 &\abs{(\xi_2,q_2)} \in I_{N_2},\abs{(\xi - \xi_1 - \xi_2,q - q_1 -q_2)} \in I_{N_3},\\
  &\abs{\tau_1 - \omega(\xi_1,q_1)} \in I_{L_1}, \abs{\tau_2 - \omega(\xi_2,q_2)} \in I_{L_2},\\ 
&\abs{\tau-\tau_1 - \tau_2 - \omega(\xi-\xi_1-\xi_2, q - q_1 -q_2)} \in I_{L_3}\}.
\end{aligned}
\end{equation}
Therefore, we have

\begin{equation}
    \label{TheNon-TrivialEstimate-WRT-A}
\Lptxy{\prod_{j=1}^3  \PQ{u}{j}}{2}{t x y}\lesssim  \sup_{(\tau,\xi,q)}\abs{A_{(\tau,\xi,q)}}^\frac{1}{2} \prod_{j=1}^3 \Lptxy{\PQ{u}{j}}{2}{t x y}.
\end{equation}
For the trivial estimate, choose the two input spatial frequencies
corresponding to $N_{\min}$ and $N_{\mathrm{med}}$; the third is determined by the
fixed output. A block
$\brs{(\xi_j,q_j):\abs{(\xi_j,q_j)}\in I_{N_j}}$ has product measure
$O(N_j^2)$. Combining this observation with the modulation-fiber argument
used below gives
\begin{equation}
\label{TheTrivialMeasureEstimate}
    \abs{A_{(\tau,\xi,q)}}
    \lesssim
    L_{\min}L_{\mathrm{med}}N_{\min}^2N_{\mathrm{med}}^2,
\end{equation}
for all $(\tau,\xi,q)\in\rrz$. To prove
\eqref{eq:trilinear-L2-N1-N3}--\eqref{eq:trilinear-L2-N1-N2},
it is enough by symmetry to prove the first estimate. Set
\begin{equation}
    \xi_3
    :=
    \xi-\xi_1-\xi_2,
    \qquad
    q_3
    :=
    q-q_1-q_2,
    \qquad
    \tau_3
    :=
    \tau-\tau_1-\tau_2.
\end{equation}
By the definition of the resonance function,
\begin{align}
&
\tau-\omega(\xi,q)
+
\mathcal H(
    \xi_1,q_1,
    \xi_2,q_2,
    \xi_3,q_3
)
\nonumber\\
&\qquad
=
\left(
    \tau_1-\omega(\xi_1,q_1)
\right)
+
\left(
    \tau_2-\omega(\xi_2,q_2)
\right)
+
\left(
    \tau_3-\omega(\xi_3,q_3)
\right).
\label{eq:resonance-modulation-identity}
\end{align}
Set $C_0:=24/5$. By the support properties of the three modulations
projections,
\begin{equation}
\left|
\tau-\omega(\xi,q)
+
\mathcal H(
    \xi_1,q_1,
    \xi_2,q_2,
    \xi_3,q_3
)
\right|
\leq
\frac85(L_1+L_2+L_3)
\leq
C_0L_{\max}.
\end{equation}
For fixed spatial frequencies, introduce
$\lambda_j:=\tau_j-\omega(\xi_j,q_j)$, $j=1,2,3$. Their sum is fixed by
\eqref{eq:resonance-modulation-identity}. The map from
$(\tau_1,\tau_2)$ to any two of the variables $\lambda_1,\lambda_2,\lambda_3$
has a Jacobian of absolute value one. Choosing the two variables corresponding
to $L_{\min}$ and $L_{\mathrm{med}}$ therefore bounds the modulation fiber by
$O(L_{\min}L_{\mathrm{med}})$. Consequently,
\begin{equation}
\label{MeasureOf-A-and-B}
    \abs{A_{(\tau,\xi,q)}}
    \lesssim
    L_{\min}L_{\mathrm{med}}\abs{B_{(\tau,\xi,q)}},
\end{equation}
where

\begin{align}
\label{Beta}
B_{(\tau,\xi,q)}
:=&
\Big\{
\big((\xi_1,q_1),(\xi_2,q_2)\big)
\in\R^2\times\Z^2:
\nonumber\\*
&
\abs{(\xi_1,q_1)}\in I_{N_1},
\quad
\abs{(\xi_2,q_2)}\in I_{N_2},
\nonumber\\*
&
\abs{
    (\xi-\xi_1-\xi_2,q-q_1-q_2)
}
\in I_{N_3},
\nonumber\\*
&
\abs{
    \tau-\omega(\xi,q)
    +
    \mathcal H(
        \xi_1,q_1,
        \xi_2,q_2,
        \xi-\xi_1-\xi_2,
        q-q_1-q_2
    )
}
\leq
C_0L_{\max}
\Big\}.
\end{align}
Recall that
\begin{equation}
    \mathcal H(
        \xi_1,q_1,
        \xi_2,q_2,
        \xi_3,q_3
    )
    =
    \omega(\xi,q)
    -
    \sum_{j=1}^3
    \omega(\xi_j,q_j),
\end{equation}
where
\begin{equation}
    \xi_3
    =
    \xi-\xi_1-\xi_2,
    \qquad
    q_3
    =
    q-q_1-q_2.
\end{equation}
Since the output $(\xi,q)$ is fixed, differentiation with respect to
$\xi_1$ gives
\begin{align}
\frac{\pt}{\pt\xi_1}
\mathcal H(
    \xi_1,q_1,
    \xi_2,q_2,
    \xi_3,q_3
)
&=
-\pt_\xi\omega(\xi_1,q_1)
+
\pt_\xi\omega(\xi_3,q_3)
\nonumber\\
&=
\abs{(\xi_3,q_3)}^2
-
\abs{(\xi_1,q_1)}^2.
\label{eq:resonance-xi1-derivative}
\end{align}
Assume that
\begin{equation}
    N_1\geq4N_3
    \qquad\text{or}\qquad
    N_3\geq4N_1.
\end{equation}
If $N_1\geq4N_3$, then the support properties of the spatial
Littlewood--Paley multipliers give
\begin{equation*}
    \abs{(\xi_1,q_1)}\geq\frac58N_1,
    \qquad
    \abs{(\xi_3,q_3)}\leq\frac85N_3\leq\frac25N_1.
\end{equation*}
Consequently,
\begin{align}
\abs{(\xi_1,q_1)}^2-\abs{(\xi_3,q_3)}^2
&\geq
\left(\frac{25}{64}-\frac4{25}\right)N_1^2
\nonumber\\
&=
\frac{369}{1600}N_1^2.
\label{eq:frequency-separation-N1-N3}
\end{align}
The same argument with the indices $1$ and $3$ interchanged applies when
$N_3\geq4N_1$. Therefore, throughout the spatial-frequency support under
consideration,
\begin{equation}
\label{eq:resonance-xi1-lower-bound}
    \left|
        \frac{\pt}{\pt\xi_1}
        \mathcal H(
            \xi_1,q_1,
            \xi_2,q_2,
            \xi_3,q_3
        )
    \right|
    \geq
    \frac{369}{1600}(N_1\vee N_3)^2.
\end{equation}
This calculation includes the borderline cases $N_1=4N_3$ and
$N_3=4N_1$.
For fixed $(q_1,q_2)\in\Z^2$, define the continuous fiber
\begin{equation}
\begin{aligned}
B_{(\tau,\xi,q)}(q_1,q_2)
:=
\Big\{
(\xi_1,\xi_2)\in\R^2:
\big((\xi_1,q_1),(\xi_2,q_2)\big)
\in B_{(\tau,\xi,q)}
\Big\}.
\end{aligned}
\end{equation}
Set $q_3:=q-q_1-q_2$. For each fixed $\xi_2\in\R$, let
\begin{equation}
\begin{aligned}
E_{\xi_2}(q_1,q_2)
:=
\Big\{
\xi_1\in\R:
&\ \abs{(\xi_1,q_1)}\in I_{N_1},
\\
&\ \abs{(\xi-\xi_1-\xi_2,q_3)}\in I_{N_3}
\Big\}.
\end{aligned}
\end{equation}
For a fixed integer $m$ and a dyadic number $M$, the set
\begin{equation*}
    \brs{r\in\R:\abs{(r,m)}\in I_M}
\end{equation*}
is the union of at most two intervals. The same remains true after a
translation or reflection of $r$. Hence $E_{\xi_2}(q_1,q_2)$ is the union
of at most four intervals.

For fixed $\xi_2$, define
\begin{equation}
\begin{aligned}
f_{\xi_2}(\xi_1)
:=
\tau-\omega(\xi,q)
+
\mathcal H(
    \xi_1,q_1,
    \xi_2,q_2,
    \xi-\xi_1-\xi_2,q_3
).
\end{aligned}
\end{equation}
On every interval component $J$ of $E_{\xi_2}(q_1,q_2)$,
\eqref{eq:resonance-xi1-lower-bound} gives
\begin{equation*}
    \inf_{\xi_1\in J}
    \abs{f_{\xi_2}'(\xi_1)}
    \geq
    \frac{369}{1600}(N_1\vee N_3)^2.
\end{equation*}
Applying \lemref{Lem:Lemma3.8} on each component with the interval
$[-C_0L_{\max},C_0L_{\max}]$ and summing over the at most four components,
we obtain
\begin{equation}
\label{eq:xi1-fiber-bound}
\abs{
\brs{
\xi_1\in E_{\xi_2}(q_1,q_2):
\abs{f_{\xi_2}(\xi_1)}\leq C_0L_{\max}
}}
\lesssim
\frac{L_{\max}}{(N_1\vee N_3)^2}.
\end{equation}
Moreover, for fixed $q_2$,
\begin{equation*}
    \abs{
    \brs{
    \xi_2\in\R:
    \abs{(\xi_2,q_2)}\in I_{N_2}
    }}
    \lesssim
    N_2.
\end{equation*}
Fubini's theorem and \eqref{eq:xi1-fiber-bound} therefore yield
\begin{equation}
\label{TheSet-B-ForFixedPoint}
    \abs{B_{(\tau,\xi,q)}(q_1,q_2)}
    \lesssim
    L_{\max}\frac{N_2}{(N_1\vee N_3)^2}.
\end{equation}
The projection of $B_{(\tau,\xi,q)}$ onto the $(q_1,q_2)$-lattice is
contained in
\begin{equation*}
| \mathcal Q_q |
:=
\left\{
(q_1,q_2)\in\Z^2:
\abs{q_1}\leq\frac85N_1,
\ \abs{q_2}\leq\frac85N_2,
\ \abs{q-q_1-q_2}\leq\frac85N_3
\right\}.
\end{equation*}
Choose the two lattice variables corresponding to the two smallest members
of $\brs{N_1,N_2,N_3}$. They have $O(N_{\min})$ and $O(N_{\mathrm{med}})$
possible values, respectively, and the remaining lattice variable is determined by
$q_1+q_2+q_3=q$. Therefore,
\begin{equation*}
    |\mathcal Q_q |
    \lesssim
    N_{\min}N_{\mathrm{med}}.
\end{equation*}
Summing \eqref{TheSet-B-ForFixedPoint} over $\mathcal Q_q$ gives
\begin{equation}
\label{TheMeasureOfSet-B}
    \abs{B_{(\tau,\xi,q)}}
    \lesssim
    L_{\max}N_{\min}N_{\mathrm{med}}
    \frac{N_2}{(N_1\vee N_3)^2}.
\end{equation}
Combining \eqref{TheMeasureOfSet-B} with \eqref{MeasureOf-A-and-B}, we obtain
\begin{equation}
\label{TheMeasureOf-A-1}
    \abs{A_{(\tau,\xi,q)}}
    \lesssim
    N_{\min}N_{\mathrm{med}}L_{\max}L_{\min}L_{\mathrm{med}}
    \frac{N_2}{(N_1\vee N_3)^2}.
\end{equation}
Since $L_{\max}L_{\min}L_{\mathrm{med}}=L_1L_2L_3$, substituting
\eqref{TheMeasureOf-A-1} into \eqref{TheNon-TrivialEstimate-WRT-A} proves
\eqref{eq:trilinear-L2-N1-N3}. The other two estimates follow
by permuting the input indices.

\end{proof}

\begin{lemma}[Trilinear estimate in Bourgain spaces]
\label{lem:trilinear-bourgain}
Let $s>1$. Then there exists $\delta>0$ sufficiently small such that
\begin{equation}
\label{eq:trilinear-bourgain}
    \xnorm{
        \pt_x\prc{\prod_{j=1}^3u_j}
    }{s}{-\frac12+3\delta}
    \lesssim
    \prod_{j=1}^3
    \xnorm{u_j}{s}{\frac12+\delta},
\end{equation}
for all $u_j\in X^{s,\frac12+\delta}$, $j=1,2,3$.
The implicit constant depends only on $s$ and $\delta$.
\end{lemma}

\begin{remark}[Comparison with previous trilinear estimates]
For $k=2$, \cite[Proposition~3.3]{farah2024note} yields the corresponding
estimate on the full timeline with output exponent
$-\frac12+2\varepsilon$ for $s>\frac56$. Its stated right-hand side places
one input at regularity $s$ and the other two at regularity $\frac56+$;
choosing $\frac56<\theta<s$ and using
$X^{s,b}\hookrightarrow X^{\theta,b}$ gives the symmetric form
\eqref{eq:farah-molinet-trilinear-estimate}. The time-localized estimate of
\cite[Proposition~4.3]{nowickikoth2026strichartz} has the same
$-\frac12+2\varepsilon$ output balance and holds for
$s>\frac{11}{24}$. Estimate \eqref{eq:trilinear-bourgain} therefore has a
higher target modulation exponent for the same input exponent, but a more
restrictive spatial regularity range. As explained in the introduction, the
exact $3\delta$ balance is not obtained from either cited estimate merely by
monotonicity of Bourgain norms.
\end{remark}

\begin{proof}
By density, it suffices to prove the estimate for Schwartz functions.
We fix throughout
\begin{equation}
\label{eq:delta-global-choice}
    0<\delta<\frac1{24}.
\end{equation}
By duality, \eqref{eq:trilinear-bourgain} follows once we prove that
\begin{equation}
\label{eq:trilinear-duality}
\begin{aligned}
    \left|
    \int_{\R\times\R\times\T}
        u_1u_2u_3\,\pt_x\overline{g}
        \,dt\,dx\,dy
    \right|
    \lesssim
    \left(
        \prod_{j=1}^3
        \xnorm{u_j}{s}{\frac12+\delta}
    \right)
    \xnorm{g}{-s}{\frac12-3\delta}
\end{aligned}
\end{equation}
for every $g\in X^{-s,\frac12-3\delta}$.

Set
\begin{equation}
\label{eq:weighted-input-profiles}
    f_j(\tau_j,\xi_j,q_j)
    :=
    \inp{\sigma_j}^{\frac12+\delta}
    \inp{\abs{(\xi_j,q_j)}}^s
    \abs{\Fh{u_j}(\tau_j,\xi_j,q_j)},
    \qquad
    j=1,2,3,
\end{equation}
and
\begin{equation}
\label{eq:weighted-output-profile}
    f_0(\tau,\xi,q)
    :=
    \inp{\sigma}^{\frac12-3\delta}
    \inp{\abs{(\xi,q)}}^{-s}
    \abs{\Fh{g}(\tau,\xi,q)}.
\end{equation}
Here
\begin{equation}
    \sigma_j
    :=
    \tau_j-\omega(\xi_j,q_j),
    \qquad
    \sigma
    :=
    \tau-\omega(\xi,q),
\end{equation}
with
\begin{equation}
    \omega(\xi,q)
    =
    \xi^3+\xi q^2.
\end{equation}
By Plancherel's identity,
\begin{equation}
\label{eq:profiles-L2-norms}
    \Lpn{f_j}{2}
    =
    \xnorm{u_j}{s}{\frac12+\delta},
    \qquad
    j=1,2,3,
\end{equation}
and
\begin{equation}
\label{eq:output-profile-L2}
    \Lpn{f_0}{2}
    =
    \xnorm{g}{-s}{\frac12-3\delta}.
\end{equation}

Let
\begin{equation}
    \tau_3=\tau-\tau_1-\tau_2,
    \qquad
    \xi_3=\xi-\xi_1-\xi_2,
    \qquad
    q_3=q-q_1-q_2,
\end{equation}
and write
\begin{equation}
    d\nu
    :=
    d\tau\,d\tau_1\,d\tau_2\,
    d\xi\,d\xi_1\,d\xi_2.
\end{equation}
The left-hand side of \eqref{eq:trilinear-duality} is bounded by
\begin{equation}
\label{eq:duality-integral-I}
    I
    :=
    \int_{\R^6}
    \sum_{q,q_1,q_2\in\Z}
    \Gamma
    \prod_{j=0}^3
    f_j(\tau_j,\xi_j,q_j)
    \,d\nu,
\end{equation}
where $(\tau_0,\xi_0,q_0)=(\tau,\xi,q)$ and
\begin{equation}
\label{TheGamma}
    \Gamma
    =
    \abs{\xi}
    \inp{\sigma}^{-\frac12+3\delta}
    \prod_{j=1}^3
    \inp{\sigma_j}^{-\frac12-\delta}
    \inp{\abs{(\xi,q)}}^s
    \prod_{j=1}^3
    \inp{\abs{(\xi_j,q_j)}}^{-s}.
\end{equation}
Thus, it remains to prove
\begin{equation}
\label{TheSufficeToProve}
    I
    \lesssim
    \prod_{j=0}^3
    \Lpn{f_j}{2}.
\end{equation}
For notational convenience, in the remainder of the proof, we denote the inverse Fourier transforms of the nonnegative profiles $f_j$
again by $u_j$. Thus
\begin{equation}
    \Lpn{u_j}{2}
    =
    \Lpn{f_j}{2},
    \qquad
    j=0,1,2,3.
\end{equation}
We use the spatial Littlewood--Paley projections $P_N$ and the
modulation projections $Q_L$, with inhomogeneous dyadic parameters
$N,L\geq1$. Write
\begin{equation}
    I
    =
    \sum_{\Ns}I_{\Ns},
\end{equation}
where
\begin{equation}
     I_{\Ns}
    =
    \sum_{\Ls}I_{\Ns}^{\Ls}.
\end{equation}
For fixed $\Ns$ and $\Ls$, set
\begin{equation}
\label{TheMain-INs}
\begin{aligned}
    J
    :=
    I_{\Ns}^{\Ls}
    =
    \int_{\R^6}
    \sum_{q,q_1,q_2}
    \Gamma
    \prod_{j=0}^3
    \abs{
        \Fh{
            P_{N_j}Q_{L_j}u_j
        }(\tau_j,\xi_j,q_j)
    }
    \,d\nu,
\end{aligned}
\end{equation}
where
\begin{equation}
    N_0=N,
    \qquad
    L_0=L.
\end{equation}
We order the input frequencies as
\begin{equation}
    N_3\leq N_2\leq N_1.
\end{equation}
For the six-case decomposition below, we use the concrete dyadic conventions
$A\ll B$ when $A<B/4$ and $A\sim B$ when $B/4\leq A\leq4B$.
Denote the largest and second-largest elements of $\{N,N_1,N_2,N_3\}$ by
$M_{\max}$ and $M_{\mathrm{next}}$, respectively. The Fourier convolution
constraint implies that $M_{\max}$ is comparable to $M_{\mathrm{next}}$.
Indeed, if the largest dyadic scale were at least eight times the second
largest, then the lower support bound $5M_{\max}/8$ for the largest frequency would exceed
the sum of the three upper support bounds, which is at most
$3(8/5)M_{\mathrm{next}}\leq3M_{\max}/5$; hence $I_{\Ns}=0$.
Because the scales are dyadic, the six cases below exhaust all remaining
alternatives.

\noindent
\textbf{\emph{Case 1: Low $\times$ Low $\times$ Low $\to$ Low interaction:}}
\begin{equation}
    N_1\leq4.
\end{equation}
By the nonvanishing condition, all four spatial frequencies are then bounded. Since
\begin{equation}
    \inp{\sigma}^{-\frac12+3\delta}
    \leq1,
\end{equation}
and all spatial Sobolev weights are uniformly bounded on this region,
we obtain
\begin{align}
I_{\Ns}
&\lesssim
\int_{\R\times\R\times\T}
\prod_{j=1}^3
\abs{
\mathcal F^{-1}
\left(
    \frac{
        \Fh{P_{N_j}u_j}
    }{
        \inp{\sigma_j}^{\frac12+\delta}
    }
\right)
}
\,
\abs{P_Nu_0}
\,dt\,dx\,dy
\nonumber\\
&\lesssim
\prod_{j=1}^3
\left\|
\mathcal F^{-1}
\left(
    \frac{
        \Fh{P_{N_j}u_j}
    }{
        \inp{\sigma_j}^{\frac12+\delta}
    }
\right)
\right\|_{L^6_{t,x,y}}
\Lpn{P_Nu_0}{2}.
\label{eq:case1-holder}
\end{align}
By Hausdorff--Young and H\"older's inequality in Fourier space,
\begin{align}
\left\|
\mathcal F^{-1}
\left(
    \frac{
        \Fh{P_{N_j}u_j}
    }{
        \inp{\sigma_j}^{\frac12+\delta}
    }
\right)
\right\|_{L^6_{t,x,y}}
\nonumber
&
\lesssim
\left\|
    \inp{\sigma_j}^{-\frac12-\delta}
    \mathbf 1_{\{\abs{(\xi_j,q_j)}\lesssim1\}}
    \Fh{P_{N_j}u_j}
\right\|_{L^{6/5}_{\tau,\xi,q}}
\nonumber\\
&
\lesssim
\left\|
    \inp{\sigma_j}^{-\frac12-\delta}
    \mathbf 1_{\{\abs{(\xi_j,q_j)}\lesssim1\}}
\right\|_{L^3_{\tau,\xi,q}}
\Lpn{\Fh{P_{N_j}u_j}}{2}
\nonumber\\
&
\lesssim
\Lpn{P_{N_j}u_j}{2}.
\label{eq:case1-L6}
\end{align}
The $L^3_{\tau,\xi,q}$ norm in the second line is finite because only finitely
many transverse modes occur and
\begin{equation}
    \int_{\R}
    \inp{\sigma}^{-\frac32-3\delta}
    \,d\sigma
    <\infty.
\end{equation}
There are only finitely many spatial dyadic scales in the present case. Hence
\begin{equation}
\label{eq:case1-final}
    \Ifr{L}{L}{L}{L}
    \lesssim
    \prod_{j=0}^3
    \Lpn{u_j}{2}.
\end{equation}

\noindent
\textbf{\emph{Case 2: High $\times$ Low $\times$ Low $\to$ High interaction:}}
\begin{equation}
    N_3\leq N_2\ll N_1\sim N.
\end{equation}
From \eqref{TheGamma}, the dyadic localization, and the previously
established estimate
\eqref{eq:trilinear-L2-N1-N2} of
\lemref{Lem:TheTrilinearEstimateIn-L2}, we obtain
\begin{align}
J
&\lesssim
L^{-\frac12+3\delta}
(L_1L_2L_3)^{-\delta}
N_2^{\frac34-s}
N_3^{\frac34-s}
\prod_{j=0}^3
\Lpn{P_{N_j}Q_{L_j}u_j}{2}.
\label{eq:case2-dyadic}
\end{align}
Since $\delta<1/6$, the modulation factors are summable. Thus
\begin{align}
\Ifr{H}{L}{L}{H}
&\lesssim
\sum_{N\sim N_1}
\Lpn{P_Nu_0}{2}
\Lpn{P_{N_1}u_1}{2}
\nonumber\\
&\quad\times
\sum_{N_2\ll N_1}
N_2^{\frac34-s}
\Lpn{P_{N_2}u_2}{2}
\left(
    \sum_{N_3\leq N_2}
    N_3^{\frac34-s}
    \Lpn{P_{N_3}u_3}{2}
\right).
\label{eq:case2-spatial-sum}
\end{align}
Since $s>1$, in particular $s>\frac34$, and therefore
\begin{equation}
    \sum_M M^{-2(s-\frac34)}
    <\infty
\end{equation}
over dyadic $M\geq1$. Cauchy--Schwarz gives
\begin{equation}
    \sum_{N_j}
    N_j^{\frac34-s}
    \Lpn{P_{N_j}u_j}{2}
    \lesssim
    \Lpn{u_j}{2},
    \qquad
    j=2,3.
\end{equation}
Finally, because $N\sim N_1$ has bounded dyadic multiplicity,
Cauchy--Schwarz in $N$ and $N_1$ yields
\begin{equation}
\label{eq:case2-final}
    \Ifr{H}{L}{L}{H}
    \lesssim
    \prod_{j=0}^3
    \Lpn{u_j}{2}.
\end{equation}

\noindent
\textbf{\emph{Case 3: High $\times$ High $\times$ Low $\to$ Low interaction:}}
\begin{equation}
    N_3,N\ll N_1\sim N_2.
\end{equation}
By \eqref{eq:trilinear-L2-N1-N3} of
\lemref{Lem:TheTrilinearEstimateIn-L2},
\begin{equation}
    \left\|
        \prod_{j=1}^3
        P_{N_j}Q_{L_j}u_j
    \right\|_2
    \lesssim
    N_3^{\frac12}
    (L_1L_2L_3)^{\frac12}
    \prod_{j=1}^3
    \Lpn{P_{N_j}Q_{L_j}u_j}{2}.
\end{equation}
Since $N_1\sim N_2$, it follows that
\begin{align}
J
&\lesssim
N_1^{-(s-1)}
N_3^{-(s-\frac12)}
L^{-\frac12+3\delta}
(L_1L_2L_3)^{-\delta}
\prod_{j=0}^3
\Lpn{P_{N_j}Q_{L_j}u_j}{2}.
\label{eq:case3-dyadic}
\end{align}
Summing the modulation parameters gives
\begin{align}
\Ifr{H}{H}{L}{L}
&\lesssim
\sum_{\substack{
N_3,N\ll N_1\\
N_1\sim N_2
}}
N_1^{-(s-1)}
N_3^{-(s-\frac12)}
\prod_{j=0}^3
\Lpn{P_{N_j}u_j}{2}.
\label{eq:case3-spatial}
\end{align}
Since $s>1$,
\begin{align}
\sum_{N_3\ll N_1}
N_3^{-(s-\frac12)}
\Lpn{P_{N_3}u_3}{2}
&\leq
\left(
    \sum_{N_3\geq1}
    N_3^{-2(s-\frac12)}
\right)^{\frac12}
\left(
    \sum_{N_3}
    \Lpn{P_{N_3}u_3}{2}^2
\right)^{\frac12}
\nonumber\\
&\lesssim
\Lpn{u_3}{2}.
\label{eq:case3-N3-sum}
\end{align}
For the low output frequency, Cauchy--Schwarz gives
\begin{align}
\sum_{N\ll N_1}
\Lpn{P_Nu_0}{2}
&\lesssim
(1+\log N_1)^{\frac12}
\Lpn{u_0}{2}.
\label{eq:case3-low-output}
\end{align}
Since $s>1$,
\begin{equation}
\label{eq:case3-log-absorption}
    \sup_{M\geq1}
    M^{-(s-1)}
    (1+\log M)^{\frac12}
    <\infty.
\end{equation}
Therefore,
\begin{align}
\Ifr{H}{H}{L}{L}
&\lesssim
\Lpn{u_0}{2}\Lpn{u_3}{2}
\sum_{N_1\sim N_2}
\Lpn{P_{N_1}u_1}{2}
\Lpn{P_{N_2}u_2}{2}
\nonumber\\
&\lesssim
\prod_{j=0}^3
\Lpn{u_j}{2}.
\label{eq:case3-final}
\end{align}
\noindent
\textbf{\emph{Case 4: High $\times$ High $\times$ Low $\to$ High interaction:}}
\begin{equation}
    N_3\ll N_1\sim N_2\sim N.
\end{equation}
Using \eqref{eq:trilinear-L2-N1-N3} of
\lemref{Lem:TheTrilinearEstimateIn-L2}, we obtain
\begin{align}
J
&\lesssim
N_1^{-(s-1)}
N_3^{-(s-\frac12)}
L^{-\frac12+3\delta}
(L_1L_2L_3)^{-\delta}
\prod_{j=0}^3
\Lpn{P_{N_j}Q_{L_j}u_j}{2}.
\label{eq:case4-dyadic}
\end{align}
After summing the modulation parameters and then $N_3$ as in
\eqref{eq:case3-N3-sum},
\begin{align}
\Ifr{H}{H}{L}{H}
&\lesssim
\Lpn{u_3}{2}
\sum_{N_1\sim N_2\sim N}
N_1^{-(s-1)}
\prod_{j=0}^2
\Lpn{P_{N_j}u_j}{2}.
\end{align}
Since $s>1$ and $N_1\geq1$,
\begin{equation}
    N_1^{-(s-1)}
    \leq1.
\end{equation}
Hence, by discrete H\"older,
\begin{align}
\Ifr{H}{H}{L}{H}
&\lesssim
\Lpn{u_3}{2}
\prod_{j=0}^2
\left(
    \sum_{N_j}
    \Lpn{P_{N_j}u_j}{2}^3
\right)^{\frac13}.
\end{align}
Since $\ell^2\subset\ell^3$,
\begin{equation}
    \left(
        \sum_M a_M^3
    \right)^{\frac13}
    \leq
    \left(
        \sum_M a_M^2
    \right)^{\frac12},
\end{equation}
and therefore
\begin{equation}
\label{eq:case4-final}
    \Ifr{H}{H}{L}{H}
    \lesssim
    \prod_{j=0}^3
    \Lpn{u_j}{2}.
\end{equation}

\noindent
\textbf{\emph{Case 5: High $\times$ High $\times$ High $\to$ Low interaction:}}
\begin{equation}
    N\ll N_1\sim N_2\sim N_3.
\end{equation}
Applying \eqref{eq:trilinear-L2-N1-N3} to the factors,
in the order $P_NQ_Lu_0$, $P_{N_2}Q_{L_2}u_2$, and
$P_{N_1}Q_{L_1}u_1$, and using $N_1 \sim N_2$, we have
\begin{equation}
\label{eq:case5-separated-bound}
\begin{aligned}
&
\left\|
    (P_NQ_Lu_0)
    (P_{N_1}Q_{L_1}u_1)
    (P_{N_2}Q_{L_2}u_2)
\right\|_2
\\
&\qquad\qquad \qquad \qquad \lesssim
N^{\frac12}
(LL_1L_2)^{\frac12}
\prod_{j=0}^2
\Lpn{P_{N_j}Q_{L_j}u_j}{2}.
\end{aligned}
\end{equation}
On the other hand, by Fourier inversion and Cauchy--Schwarz,
\begin{equation}
    \norm{P_{N_j}Q_{L_j}u_j}_{L^\infty}
    \lesssim
    N_jL_j^{\frac12}
    \norm{P_{N_j}Q_{L_j}u_j}_{L^2},
    \qquad j=1,2.
\end{equation}
Hence, since $N_1\sim N_2$,
\begin{equation}
\begin{aligned}
\norm{
    (P_NQ_Lu_0)
    (P_{N_1}Q_{L_1}u_1)
    (P_{N_2}Q_{L_2}u_2)
}_{L^2}
\lesssim
N_1^2L_1^{\frac12}L_2^{\frac12}
\prod_{j=0}^2
\norm{P_{N_j}Q_{L_j}u_j}_{L^2}.
\end{aligned}
\end{equation}
Taking the geometric mean of the preceding estimate and
\eqref{eq:case5-separated-bound}, with weights $\theta$ and
$1-\theta$, respectively, gives

\begin{equation}
\label{eq:case5-interpolated}
\begin{aligned}
&
\left\|
    (P_NQ_Lu_0)
    (P_{N_1}Q_{L_1}u_1)
    (P_{N_2}Q_{L_2}u_2)
\right\|_2
\\
&\qquad\qquad \qquad \qquad \lesssim
N^{\frac{1-\theta}{2}}
N_1^{2\theta}
L^{\frac{1-\theta}{2}}
L_1^{\frac12}
L_2^{\frac12}
\prod_{j=0}^2
\Lpn{P_{N_j}Q_{L_j}u_j}{2}.
\end{aligned}
\end{equation}
Fix $\theta=\frac14$. Then, using
\begin{equation}
    N\ll N_1\sim N_2\sim N_3
\end{equation}
and the spatial factor in \eqref{TheGamma}, we deduce
\begin{align}
J
&\lesssim
N_1^{-\alpha}
L^{3\delta-\frac18}
L_1^{-\delta}
L_2^{-\delta}
L_3^{-\frac12-\delta}
\prod_{j=0}^3
\Lpn{P_{N_j}Q_{L_j}u_j}{2},
\label{eq:case5-dyadic}
\end{align}
where $ \alpha:=2s-\frac{15}{8}$. Since $s>1$,
\begin{equation}
    \alpha>\frac18>0.
\end{equation}
Moreover, by \eqref{eq:delta-global-choice},
\begin{equation}
    3\delta-\frac18<0,
\end{equation}
and therefore all four modulation weights in
\eqref{eq:case5-dyadic} are summable.

Consequently,
\begin{align}
\Ifr{H}{H}{H}{L}
&\lesssim
\sum_{\substack{
N\ll N_1\\
N_1\sim N_2\sim N_3
}}
N_1^{-\alpha}
\prod_{j=0}^3
\Lpn{P_{N_j}u_j}{2}.
\label{eq:case5-spatial}
\end{align}
For fixed $N_1$,
\begin{equation}
    \sum_{N\ll N_1}
    \Lpn{P_Nu_0}{2}
    \lesssim
    (1+\log N_1)^{\frac12}
    \Lpn{u_0}{2}.
\end{equation}
Since $\alpha>0$,
\begin{equation}
    \sup_{M\geq1}
    M^{-\alpha}
    (1+\log M)^{\frac12}
    <\infty.
\end{equation}
Thus
\begin{align}
\Ifr{H}{H}{H}{L}
&\lesssim
\Lpn{u_0}{2}
\sum_{N_1\sim N_2\sim N_3}
\prod_{j=1}^3
\Lpn{P_{N_j}u_j}{2}
\nonumber\\
&\lesssim
\Lpn{u_0}{2}
\prod_{j=1}^3
\left(
    \sum_{N_j}
    \Lpn{P_{N_j}u_j}{2}^3
\right)^{\frac13}
\nonumber\\
&\lesssim
\prod_{j=0}^3
\Lpn{u_j}{2}.
\label{eq:case5-final}
\end{align}
\medskip
\noindent
\textbf{\emph{Case 6: High $\times$ High $\times$ High $\to$ High interaction:}}
\begin{equation}
\label{eq:case6-frequency-comparability}
    N_1\geq8,
    \qquad
    N_1\sim N_2\sim N_3\sim N.
\end{equation}
For fixed $N,N_1,N_2,N_3,L,L_1,L_2,L_3$, define
\begin{equation}
\label{eq:case6-Fj}
    F_j
    :=
    \abs{
        \Fh{
            P_{N_j}Q_{L_j}u_j
        }
    },
    \qquad
    v_j
    :=
    \mathcal F^{-1}F_j,
    \qquad
    j=0,1,2,3.
\end{equation}
By Plancherel's identity,
\begin{equation}
\label{eq:case6-majorant-L2}
    \Lpn{v_j}{2}
    =
    \Lpn{P_{N_j}Q_{L_j}u_j}{2}.
\end{equation}
The replacement of each Fourier transform by its absolute value is legitimate because the Fourier integrand defining $I$ is nonnegative.  We divide the argument into the regions
\begin{equation}
    \abs{\xi}\leq100
\end{equation}
and
\begin{equation}
    \abs{\xi}>100.
\end{equation}
\medskip
\noindent
\textbf{\emph{The region $\abs{\xi}\leq100$:}}
For $k\geq0$, let
\begin{equation}
    E_k
    :=
    \brs{
        (100)\,2^{-(k+1)}
        \leq
        \abs{\xi}
        \leq
        (100)\,2^{-k}
    }.
\end{equation}
The set $\{\xi=0\}$ is negligible. Let $J_k$ denote the contribution of $E_k$ to $J$. On $E_k$, we have
\begin{align}
\abs{\xi} \inp{\abs{(\xi,q)}}^s 
    \displaystyle\prod_{j=1}^3
    \inp{\abs{(\xi_j,q_j)}}^{-s}
\lesssim
2^{-k}
\frac{N^s}{N_1^sN_2^sN_3^s}
\lesssim
2^{-k}N^{-2s}.
\label{eq:case6-small-xi-spatial-multiplier}
\end{align}
After extracting this factor, we may enlarge the region of
integration by discarding the condition $\xi\in E_k$.
Hence, by Plancherel and Cauchy--Schwarz,
\begin{align}
J_k
&\lesssim
2^{-k}N^{-2s}
L^{-\frac12+3\delta}
\prod_{j=1}^3L_j^{-\frac12-\delta}
\Lpn{v_1v_2v_3}{2}
\Lpn{v_0}{2}.
\label{eq:case6-small-xi-before-trivial}
\end{align}
For later use, we record the elementary estimate
\begin{equation}
\label{eq:case6-trivial-Linfty}
    \Lpn{v_j}{\infty}
    \lesssim
    N_jL_j^{\frac12}
    \Lpn{v_j}{2}.
\end{equation}
Indeed,
\begin{equation}
    \abs{
        \operatorname{supp}\Fh{v_j}
    }
    \lesssim
    N_j^2L_j,
\end{equation}
because the $\xi_j$-interval has length $O(N_j)$, there are
$O(N_j)$ admissible integers $q_j$, and the $\tau_j$-fiber has
length $O(L_j)$. Thus, Fourier inversion and Cauchy--Schwarz give
\eqref{eq:case6-trivial-Linfty}. Using \eqref{eq:case6-frequency-comparability},
\begin{align}
\Lpn{v_1v_2v_3}{2}
&\leq
\Lpn{v_1}{\infty}
\Lpn{v_2}{\infty}
\Lpn{v_3}{2}
\nonumber\\
&\lesssim
N^2
(L_1L_2)^{\frac12}
\prod_{j=1}^3
\Lpn{v_j}{2}.
\label{eq:case6-trivial-trilinear}
\end{align}
Substitution into \eqref{eq:case6-small-xi-before-trivial} gives
\begin{align}
J_k
&\lesssim
2^{-k}
N^{2-2s}
L^{-\frac12+3\delta}
(L_1L_2)^{-\delta}
L_3^{-\frac12-\delta}
\prod_{j=0}^3
\Lpn{P_{N_j}Q_{L_j}u_j}{2}.
\label{eq:case6-small-xi-dyadic}
\end{align}
Since $\delta<1/24<1/6$, all modulation weights are summable, and
\begin{equation}
    \sum_{k\geq0}2^{-k}<\infty.
\end{equation}
Therefore,
\begin{equation}
\label{eq:case6-small-xi-summed}
    \sum_{k\geq0}
    \sum_{\Ls}J_k
    \lesssim
    N^{2-2s}
    \prod_{j=0}^3
    \Lpn{P_{N_j}u_j}{2}.
\end{equation}
Since $s>1$,
\begin{equation}
    N^{2-2s}\leq1.
\end{equation}
The bounded multiplicity of
$N_1\sim N_2\sim N_3\sim N$, followed by discrete H\"older and
$\ell^2\subset\ell^4$, yields
\begin{equation}
\label{eq:case6-small-xi-final}
    \sum_{N_1\sim N_2\sim N_3\sim N}
    \sum_{k\geq0}
    \sum_{\Ls}
    J_k
    \lesssim
    \prod_{j=0}^3
    \Lpn{u_j}{2}.
\end{equation}

\medskip
\noindent
\emph{A pair-frequency-localized bilinear estimate.}
The mechanism and the $K^{-\frac14}$ gain below are obtained by retaining the
longitudinal output scale in the high-output-frequency branch of the
quadratic-phase argument in \cite[proof of Proposition~3.1]{farah2024note}.
We record the estimate in this pair-frequency-localized form because the
argument below keeps the parameter $K$ and applies the bound to two paired
convolutions. For a dyadic number $K\geq2$, define
\begin{equation}
    \mathcal A_K
    :=
    \brs{
        \eta\in\R:
        K\leq\abs{\eta}<2K
    }.
\end{equation}
Choose once and for all an even function $\chi\in C_0^\infty(\R)$ such that
$0\leq\chi\leq1$, $\chi=1$ on $\{1\leq\abs r\leq2\}$, and
$\supp(\chi)\subset\{\frac12\leq\abs r\leq4\}$. Set
\begin{equation}
    \chi_K(\eta):=\chi(\eta/K).
\end{equation}
Define $R_K$ by
\begin{equation}
    \mathcal F(R_Kw)(\rho,\eta,m)
    :=
    \chi_K(\eta)\Fh{w}(\rho,\eta,m).
\end{equation}
We claim that, if $\Fh{w_j}$ is supported in
\begin{equation}
    \abs{(\eta_j,m_j)}\sim M_j,
    \qquad
    \abs{
        \rho_j-\omega(\eta_j,m_j)
    }
    \lesssim H_j,
    \qquad
    j=1,2,
\end{equation}
then
\begin{equation}
\label{eq:case6-bilinear-estimate}
    \boxed{
    \Lpn{R_K(w_1w_2)}{2}
    \lesssim
    M_{\min}^{\frac12}
    K^{-\frac14}
    H_{\min}^{\frac12}
    H_{\max}^{\frac14}
    \prod_{j=1}^2
    \Lpn{w_j}{2}
    },
\end{equation}
where
\begin{equation}
    M_{\min}:=M_1\wedge M_2,
    \qquad
    H_{\min}:=H_1\wedge H_2,
    \qquad
    H_{\max}:=H_1\vee H_2.
\end{equation}
Fix $(\rho,\eta,m)$ satisfying
\begin{equation}
    \frac K2
    \leq
    \abs{\eta}
    \leq
    4K.
\end{equation}
Let $A_{\rho,\eta,m}$ be the set of
$(\rho_1,\eta_1,m_1)\in\R^2\times\Z$ such that
\begin{equation}
\begin{aligned}
    &\abs{(\eta_1,m_1)}\sim M_1,
    \qquad
    \abs{(\eta-\eta_1,m-m_1)}\sim M_2,
    \\
    &\abs{
        \rho_1-\omega(\eta_1,m_1)
    }
    \lesssim H_1,
    \\
    &\abs{
        \rho-\rho_1
        -
        \omega(\eta-\eta_1,m-m_1)
    }
    \lesssim H_2.
\end{aligned}
\end{equation}
The convolution Cauchy--Schwarz inequality and Plancherel give
\begin{equation}
\label{eq:case6-convolution-CS}
    \Lpn{R_K(w_1w_2)}{2}
    \lesssim
    \sup_{\substack{
        (\rho,\eta,m)\\
        K/2\leq\abs{\eta}\leq4K
    }}
    \abs{A_{\rho,\eta,m}}^{\frac12}
    \prod_{j=1}^2
    \Lpn{w_j}{2}.
\end{equation}
For fixed $(\eta_1,m_1)$, the admissible values of $\rho_1$ must satisfy
\begin{equation}
    \abs{\rho_1-\omega(\eta_1,m_1)}
    \lesssim H_1
\end{equation}
and
\begin{equation}
    \abs{
        \rho-\rho_1-\omega(\eta-\eta_1,m-m_1)
    }
    \lesssim H_2.
\end{equation}
Thus, $\rho_1$ belongs to the intersection of two intervals of lengths
$O(H_1)$ and $O(H_2)$, respectively. Consequently,
\begin{equation}
\label{eq:case6-rho1-fiber}
\left|
\left\{
\rho_1\in\R:
\begin{array}{l}
\abs{\rho_1-\omega(\eta_1,m_1)}\lesssim H_1,\\[1mm]
\abs{\rho-\rho_1-\omega(\eta-\eta_1,m-m_1)}\lesssim H_2
\end{array}
\right\}
\right|
\lesssim
H_{\min}.
\end{equation}
Moreover, the two modulation restrictions imply
\begin{equation}
\label{eq:case6-resonance-restriction}
    \abs{
        \rho
        -
        \omega(\eta_1,m_1)
        -
        \omega(\eta-\eta_1,m-m_1)
    }
    \lesssim H_{\max}.
\end{equation}
For fixed $m_1$, set
\begin{equation}
    h(\eta_1)
    :=
    \omega(\eta_1,m_1)
    +
    \omega(\eta-\eta_1,m-m_1).
\end{equation}
Since
\begin{equation}
    \omega(\eta,m)
    =
    \eta^3+\eta m^2,
\end{equation}
we have
\begin{align}
h(\eta_1)
&=
3\eta\eta_1^2
+
\left(
    m_1^2
    -(m-m_1)^2
    -3\eta^2
\right)\eta_1
\nonumber\\
&\qquad
+
\eta^3
+
\eta(m-m_1)^2,
\end{align}
and hence
\begin{equation}
\label{eq:case6-second-derivative}
    h''(\eta_1)=6\eta.
\end{equation}
Thus
\begin{equation}
    \abs{h''(\eta_1)}
    \sim K.
\end{equation}
For a quadratic polynomial $p(x)=ax^2+bx+c$ and an interval
$\mathcal I\subset\R$,
\begin{equation}
\label{eq:case6-quadratic-sublevel}
    \abs{
        \brs{
            x:p(x)\in\mathcal I
        }
    }
    \lesssim
    \left(
        \frac{\abs{\mathcal I}}{\abs a}
    \right)^{\frac12}.
\end{equation}
This follows immediately by completing the square and reducing to
the estimate
\begin{equation}
    \abs{
        \brs{
            x:x^2\in[\alpha,\beta]
        }
    }
    \lesssim
    \sqrt{\beta-\alpha}.
\end{equation}
Applying \eqref{eq:case6-quadratic-sublevel} to
\eqref{eq:case6-resonance-restriction}, we obtain
\begin{equation}
\label{eq:case6-eta1-fiber}
    \abs{
        \brs{
            \eta_1:
            \eqref{eq:case6-resonance-restriction}
            \text{ holds}
        }
    }
    \lesssim
    \left(
        \frac{H_{\max}}{K}
    \right)^{\frac12}.
\end{equation}
Finally,
\begin{equation}
    \abs{m_1}\lesssim M_1,
    \qquad
    \abs{m-m_1}\lesssim M_2,
\end{equation}
so the number of admissible integers $m_1$ is bounded by
\begin{equation}
\label{eq:case6-discrete-count}
    |\left\{m_1\in\Z:
        \abs{m_1}\lesssim M_1,
        \ \abs{m-m_1}\lesssim M_2
    \right\}|
    \lesssim
    M_{\min}.
\end{equation}
Combining
\eqref{eq:case6-rho1-fiber},
\eqref{eq:case6-eta1-fiber}, and
\eqref{eq:case6-discrete-count}, we get
\begin{equation}
    \abs{A_{\rho,\eta,m}}
    \lesssim
    M_{\min}
    H_{\min}
    \left(
        \frac{H_{\max}}{K}
    \right)^{\frac12}.
\end{equation}
Substitution into \eqref{eq:case6-convolution-CS} proves
\eqref{eq:case6-bilinear-estimate}.

\medskip
\noindent
\textbf{\emph{The region $\abs{\xi}>100$:}}

Define
\begin{equation}
    \eta_{12}:=\xi_1+\xi_2,
    \qquad
    \eta_{13}:=\xi_1+\xi_3,
    \qquad
    \eta_{23}:=\xi_2+\xi_3.
\end{equation}
Since
\begin{equation}
    2\xi
    =
    \eta_{12}+\eta_{13}+\eta_{23},
\end{equation}
we have
\begin{equation}
\label{eq:case6-dominant-pair-lower-bound}
    \max
    \brs{
        \abs{\eta_{12}},
        \abs{\eta_{13}},
        \abs{\eta_{23}}
    }
    \geq
    \frac23\abs{\xi}.
\end{equation}
Partition the region into three measurable subsets according to which
pair realizes this maximum, resolving ties in a fixed order.
By symmetry in the three inputs, it suffices to consider
\begin{equation}
    \abs{\eta_{12}}
    =
    \max
    \brs{
        \abs{\eta_{12}},
        \abs{\eta_{13}},
        \abs{\eta_{23}}
    }.
\end{equation}
Decompose further into dyadic shells
\begin{equation}
    \eta_{12}\in\mathcal A_K.
\end{equation}
Then
\begin{equation}
\label{eq:case6-K-relations}
    \abs{\xi}\lesssim K,
    \qquad
    1\ll K\lesssim N.
\end{equation}
Let $J_K$ denote the corresponding contribution. The spatial part of
the multiplier satisfies
\begin{equation}
\abs{\xi} \inp{\abs{(\xi,q)}}^s \prod_{j=1}^3
    \inp{\abs{(\xi_j,q_j)}}^{-s}
\lesssim
K
\frac{N^s}{N_1^sN_2^sN_3^s}
\lesssim
KN^{-2s}.
\label{eq:case6-large-xi-spatial}
\end{equation}
Define
\begin{equation}
    F_0^\sharp(\tau,\xi,q)
    :=
    F_0(-\tau,-\xi,-q).
\end{equation}
After extracting \eqref{eq:case6-large-xi-spatial}, we may enlarge
the integration domain by discarding the condition that
$\abs{\eta_{12}}$ is the largest pair frequency, while retaining
$\eta_{12}\in\mathcal A_K$. The resulting convolution integral is
bounded by
\begin{equation}
J_K
\lesssim
KN^{-2s}
L^{-\frac12+3\delta}
\prod_{j=1}^3
L_j^{-\frac12-\delta}
\,
\norm{
    \chi_K(\eta)(F_1*F_2)
}_{L^2}
\norm{
    \chi_K(\eta)(F_0^\sharp*F_3)
}_{L^2}.
\label{eq:case6-paired-bound}
\end{equation}
Indeed, the first convolution has frequency
\begin{equation}
    (\tau_1+\tau_2,\xi_1+\xi_2,q_1+q_2),
\end{equation}
whereas
\begin{equation}
\begin{aligned}
(-\tau,-\xi,-q)
+
(\tau_3,\xi_3,q_3)
=
-(\tau_1+\tau_2,\xi_1+\xi_2,q_1+q_2).
\end{aligned}
\end{equation}
Since
\begin{equation}
    \omega(-\xi,-q)
    =
    -\omega(\xi,q),
\end{equation}
the reflected factor $F_0^\sharp$ has the same spatial-frequency and
modulation localizations as $F_0$, and
\begin{equation}
    \norm{F_0^\sharp}_{L^2}
    =
    \norm{F_0}_{L^2}.
\end{equation}
By Plancherel, apply \eqref{eq:case6-bilinear-estimate} first to functions
whose Fourier transforms are $F_1$ and $F_2$, and then to functions whose
Fourier transforms are $F_0^\sharp$ and $F_3$. We obtain
\begin{align}
\norm{
    \chi_K(\eta)(F_1*F_2)
}_{L^2}
&\lesssim
N^{\frac12}
K^{-\frac14}
(L_1\wedge L_2)^{\frac12}
(L_1\vee L_2)^{\frac14}
\prod_{j=1}^2
\Lpn{P_{N_j}Q_{L_j}u_j}{2},
\label{eq:case6-first-pair}
\end{align}
and
\begin{align}
\norm{
    \chi_K(\eta)(F_0^\sharp*F_3)
}_{L^2}
&\lesssim
N^{\frac12}
K^{-\frac14}
(L\wedge L_3)^{\frac12}
(L\vee L_3)^{\frac14}
\nonumber\\
&\qquad\times
\Lpn{P_NQ_Lu_0}{2}
\Lpn{P_{N_3}Q_{L_3}u_3}{2}.
\label{eq:case6-second-pair}
\end{align}
Consequently,
\begin{equation}
\label{eq:case6-after-bilinear}
    J_K
    \lesssim
    K^{\frac12}
    N^{1-2s}
    \mathcal M(L,L_1,L_2,L_3)
    \prod_{j=0}^3
    \Lpn{P_{N_j}Q_{L_j}u_j}{2},
\end{equation}
where
\begin{align}
\mathcal M(L,L_1,L_2,L_3)
&:=
L^{-\frac12+3\delta}
\prod_{j=1}^3
L_j^{-\frac12-\delta}
\nonumber\\
&\qquad\times
(L_1\wedge L_2)^{\frac12}
(L_1\vee L_2)^{\frac14}
\nonumber\\
&\qquad\times
(L\wedge L_3)^{\frac12}
(L\vee L_3)^{\frac14}.
\label{eq:case6-M-factor}
\end{align}
We claim that
\begin{equation}
\label{eq:case6-M-summable}
    \sum_{L,L_1,L_2,L_3}
    \mathcal M(L,L_1,L_2,L_3)
    <\infty.
\end{equation}
For the $(L_1,L_2)$ pair, if $L_1\leq L_2$, then
\begin{equation}
    L_1^{-\frac12-\delta}
    L_2^{-\frac12-\delta}
    L_1^{\frac12}
    L_2^{\frac14}
    =
    L_1^{-\delta}
    L_2^{-\frac14-\delta},
\end{equation}
which is dyadically summable; the region $L_2\leq L_1$ is
identical.

For the $(L,L_3)$ pair, if $L\leq L_3$, then
\begin{equation}
    L^{-\frac12+3\delta}
    L_3^{-\frac12-\delta}
    L^{\frac12}
    L_3^{\frac14}
    =
    L^{3\delta}
    L_3^{-\frac14-\delta}.
\end{equation}
Hence
\begin{align}
\sum_{L_3}
\sum_{L\leq L_3}
L^{3\delta}
L_3^{-\frac14-\delta}
&\lesssim
\sum_{L_3}
L_3^{-\frac14+2\delta}
<\infty,
\end{align}
because $\delta<1/24<1/8$.

If $L_3\leq L$, then
\begin{equation}
    L^{-\frac12+3\delta}
    L_3^{-\frac12-\delta}
    L_3^{\frac12}
    L^{\frac14}
    =
    L^{-\frac14+3\delta}
    L_3^{-\delta}.
\end{equation}
Since
\begin{equation}
    \sum_{L_3\leq L}
    L_3^{-\delta}
    \lesssim1,
\end{equation}
we have
\begin{equation}
    \sum_L
    L^{-\frac14+3\delta}
    <\infty,
\end{equation}
because $\delta<1/24<1/12$. This proves
\eqref{eq:case6-M-summable}.

Therefore,
\begin{equation}
\label{eq:case6-after-modulation}
    \sum_{\Ls}J_K
    \lesssim
    K^{\frac12}
    N^{1-2s}
    \prod_{j=0}^3
    \Lpn{P_{N_j}u_j}{2}.
\end{equation}
Using $K\lesssim N$,
\begin{equation}
    \sum_{K\lesssim N}
    K^{\frac12}
    \lesssim
    N^{\frac12},
\end{equation}
and hence
\begin{equation}
\label{eq:case6-after-K}
    \sum_{K\lesssim N}
    \sum_{\Ls}J_K
    \lesssim
    N^{\frac32-2s}
    \prod_{j=0}^3
    \Lpn{P_{N_j}u_j}{2}.
\end{equation}
Since $s>1$,
\begin{equation}
    \frac32-2s<-\frac12,
\end{equation}
so
\begin{equation}
    N^{\frac32-2s}
    \lesssim1.
\end{equation}
Finally, by bounded dyadic multiplicity, discrete H\"older, and
$\ell^2\subset\ell^4$,
\begin{align}
\sum_{N_1\sim N_2\sim N_3\sim N}
\prod_{j=0}^3
\Lpn{P_{N_j}u_j}{2}
&\lesssim
\prod_{j=0}^3
\left(
    \sum_{N_j}
    \Lpn{P_{N_j}u_j}{2}^4
\right)^{\frac14}
\nonumber\\
&\lesssim
\prod_{j=0}^3
\Lpn{u_j}{2}.
\end{align}
The regions in which $\abs{\eta_{13}}$ or $\abs{\eta_{23}}$ is
maximal are identical after permuting the three input indices.
Thus
\begin{equation}
\label{eq:case6-final}
    \Ifr{H}{H}{H}{H}
    \lesssim
    \prod_{j=0}^3
    \Lpn{u_j}{2}.
\end{equation}
Combining Cases $1$--$6$ proves
\begin{equation}
    I
    \lesssim
    \prod_{j=0}^3
    \Lpn{u_j}{2}.
\end{equation}
Returning to the weighted profiles and using
\eqref{eq:profiles-L2-norms}--\eqref{eq:output-profile-L2}, we obtain
\begin{equation}
\begin{aligned}
    \left|
    \int_{\R\times\R\times\T}
        u_1u_2u_3\,\pt_x\overline g
        \,dt\,dx\,dy
    \right|
    &\lesssim
    \left(
        \prod_{j=1}^3
        \xnorm{u_j}{s}{\frac12+\delta}
    \right)
    \xnorm{g}{-s}{\frac12-3\delta}.
\end{aligned}
\end{equation}
By duality, this is precisely
\eqref{eq:trilinear-bourgain}, and the proof is complete.
\end{proof}

\section{Proof of Theorem~\ref{MainResult}}
\begin{proof}
Let $s>1$, $0<T<1$, and choose $0<\delta<\frac1{24}$ as in
\lemref{lem:trilinear-bourgain}. Fix a real-valued function
$\eta\in C_0^\infty(\R)$ such that $\eta=1$ on $[-1,1]$ and
$\supp(\eta)\subset[-2,2]$.

We first record the continuous embedding of the restricted Bourgain space.
Let $u\in X_T^{s,\frac12+\delta}$ and $\varepsilon>0$. By definition of the
restricted norm, there is an extension $U\in X^{s,\frac12+\delta}$ satisfying
\begin{equation}
    \xnorm{U}{s}{\frac12+\delta}
    \leq
    \xtnorm{u}{s}{\frac12+\delta}+\varepsilon.
\end{equation}
Since
\begin{equation}
    \xnorm{U}{s}{\frac12+\delta}
    =
    \norm{
        e^{t\pt_x\Delta}U
    }_{H_t^{\frac12+\delta}H_{xy}^s},
\end{equation}
the Sobolev embedding
\begin{equation}
    H^{\frac12+\delta}(\R)
    \hookrightarrow
    C(\R)
\end{equation}
implies
\begin{equation}
    \sup_{t\in\R}
    \norm{
        e^{t\pt_x\Delta}U(t)
    }_{H^s}
    \lesssim
    \xnorm{U}{s}{\frac12+\delta}.
\end{equation}
The group $e^{-t\pt_x\Delta}$ is unitary on $H^s(\rt)$. Restricting to
$[0,T]$ and letting $\varepsilon\downarrow0$ therefore gives
\begin{equation}
    X_T^{s,\frac12+\delta}
    \hookrightarrow
    C\prc{[0,T];H^s(\rt)}.
\end{equation}

Define the integral-equation map
\begin{equation*}
    \Phi:
    H^s(\rt)\times X_T^{s,\frac12+\delta}
    \longrightarrow
    X_T^{s,\frac12+\delta}
\end{equation*}
by
\begin{equation}
\label{IntegralEquation}
\Phi(u_0,u)
:=
\eta(t)e^{-t\pt_x\Delta}u_0
-
\eta\prc{\frac{t}{T}}
\int_0^t
e^{-(t-t')\pt_x\Delta}
\pt_x\prc{u^3(t')}
\,dt'.
\end{equation}
Because both cutoffs equal one on $[0,T]$, the fixed points of $\Phi$ are
exactly the solutions of \eqref{eq:mZK-cylinder-intro} on that interval.

For $r,R>0$, set
\begin{equation*}
\begin{aligned}
    B_r^H
    &:=\{u_0\in H^s(\rt):\norm{u_0}_{H^s}<r\},
    \\
    \overline B_r^H
    &:=\{u_0\in H^s(\rt):\norm{u_0}_{H^s}\leq r\},
    \\
    \overline B_R^X
    &:=\{u\in X_T^{s,\frac12+\delta}:
        \xtnorm{u}{s}{\frac12+\delta}\leq R\}.
\end{aligned}
\end{equation*}
Let $u_0\in\overline B_r^H$ and $u\in\overline B_R^X$. For every
$\varepsilon>0$, choose an extension $w\in X^{s,\frac12+\delta}$ such that
\begin{equation}
    w|_{(0,T)\times\rt}=u
\end{equation}
and
\begin{equation}
    \norm{w}_{X^{s,\frac12+\delta}}
    \leq
    \norm{u}_{X_T^{s,\frac12+\delta}}
    +\varepsilon.
\end{equation}
Let $\widetilde\Phi(u_0,w)$ denote the full-line expression on the right-hand
side of \eqref{IntegralEquation}, with $w$ in place of $u$. Then
$\widetilde\Phi(u_0,w)|_{(0,T)\times\rt}=\Phi(u_0,u)$, and
\begin{align*}
     \xnorm{\widetilde\Phi(u_0,w)}{s}{\frac{1}{2}+\delta}& \leq  \xnorm{ \eta(t)     e^{- t \pt_x \Delta} u_0 }{s}{\frac{1}{2}+\delta}\\
     &\qquad \qquad +  \xnorm{   \eta(\frac{t}{T}) \int_0^t e^{- (t-t') \pt_x \Delta} \pt_x (w^3 (t')) dt'}{s}{\frac{1}{2}+\delta}.
\end{align*}
Here $b=\frac12+\delta$ and $b'=-\frac12+3\delta$ satisfy
$-\frac12<b'\leq0\leq b\leq b'+1$, and
$1-b+b'=2\delta$. Thus
\lemref{Lem:homogeneous-linear-estimate},
\lemref{lem:duhamel-linear-estimate}, and
\lemref{lem:trilinear-bourgain} give
\begin{align*}
      \xnorm{\widetilde\Phi(u_0,w)}{s}{\frac12+\delta}
      &\leq c_1\Hnxy{u_0}{s}{}
      +c_2c_3T^{2\delta}\xnorm{w}{s}{\frac12+\delta}^3.
\end{align*}
The definition of the restricted norm, therefore, gives
\begin{align*}
\xtnorm{\Phi(u_0,u)}{s}{\frac12+\delta}
&\leq
\xnorm{\widetilde\Phi(u_0,w)}{s}{\frac12+\delta}
\\
&\leq
c_1\Hnxy{u_0}{s}{}
+
c_2c_3T^{2\delta}
\left(
    \xtnorm{u}{s}{\frac12+\delta}
    +\varepsilon
\right)^3.
\end{align*}
Let $c\geq\max\{c_1,c_2c_3\}$ be large enough to dominate the constants in the difference estimate below. Letting $\varepsilon\downarrow0$, we obtain
\begin{equation}
    \begin{aligned}
\xtnorm{\Phi(u_0,u)}{s}{\frac12+\delta}
&\leq
c_1\Hnxy{u_0}{s}{}
+
c_2c_3T^{2\delta}
\xtnorm{u}{s}{\frac12+\delta}^3\\
&\leq cr+cT^{2\delta}R^3.
\end{aligned}
\end{equation}
Put $R=2cr$ and require
\begin{equation}
    T^{2\delta}<\frac1{24c^3r^2}.
\end{equation}
Then
\begin{equation*}
    cr+cT^{2\delta}R^3
    <\frac R2+\frac R6<R,
\end{equation*}
so
$\Phi(u_0,\cdot)$ maps $\overline B_R^X$ into itself, uniformly for
$u_0\in\overline B_r^H$.

To prove contraction, let $f_1, g_1 \in \overline B_R^X$. For every
$\varepsilon>0$, choose extensions $f,g,h\in X^{s,\frac12+\delta}$ such that
\begin{equation}
    f|_{(0,T)\times\rt}=f_1,
    \qquad
    g|_{(0,T)\times\rt}=g_1,
    \qquad
    h|_{(0,T)\times\rt}=f_1-g_1,
\end{equation}
and
\begin{equation}
\begin{aligned}
    \xnorm{f}{s}{\frac12+\delta}
    &\leq
    \xtnorm{f_1}{s}{\frac12+\delta}
    +\varepsilon,
    \\
    \xnorm{g}{s}{\frac12+\delta}
    &\leq
    \xtnorm{g_1}{s}{\frac12+\delta}
    +\varepsilon,
    \\
    \xnorm{h}{s}{\frac12+\delta}
    &\leq
    \xtnorm{f_1-g_1}{s}{\frac12+\delta}
    +\varepsilon.
\end{aligned}
\end{equation}
On $(0,T)\times\rt$,
\begin{equation}
    f_1^3-g_1^3
    =
    (f_1-g_1)(f_1^2+f_1g_1+g_1^2),
\end{equation}
and therefore
\begin{equation}
    h(f^2+fg+g^2)
\end{equation}
is an extension of $f_1^3-g_1^3$.

Using
\lemref{lem:duhamel-linear-estimate}
and
\lemref{lem:trilinear-bourgain},
we obtain
\begin{align*}
&
\xtnorm{
    \Phi(u_0,f_1)-\Phi(u_0,g_1)
}{s}{\frac12+\delta}
\\
&\qquad\lesssim
T^{2\delta}
\Big(
    \xnorm{f}{s}{\frac12+\delta}^2
    +
    \xnorm{f}{s}{\frac12+\delta}
    \xnorm{g}{s}{\frac12+\delta}
    +
    \xnorm{g}{s}{\frac12+\delta}^2
\Big)
\xnorm{h}{s}{\frac12+\delta}.
\end{align*}
Letting $\varepsilon\downarrow0$, we obtain
\begin{align*}
&
\xtnorm{
    \Phi(u_0,f_1)-\Phi(u_0,g_1)
}{s}{\frac12+\delta}
\\
&\qquad\lesssim
T^{2\delta}
\Big(
    \xtnorm{f_1}{s}{\frac12+\delta}^2
    +
    \xtnorm{f_1}{s}{\frac12+\delta}
    \xtnorm{g_1}{s}{\frac12+\delta}
    +
    \xtnorm{g_1}{s}{\frac12+\delta}^2
\Big)
\xtnorm{f_1-g_1}{s}{\frac12+\delta}
\\
&\qquad\leq
3cT^{2\delta}R^2
\xtnorm{f_1-g_1}{s}{\frac12+\delta}
\\
&\qquad<
\frac12
\xtnorm{f_1-g_1}{s}{\frac12+\delta}.
\end{align*}
Thus $\Phi(u_0,\cdot)$ is a strict contraction on $\overline B_R^X$,
uniformly for $u_0\in\overline B_r^H$. The contraction mapping theorem
defines the data-to-solution map
\begin{equation*}
    \mathcal S_r:B_r^H\longrightarrow\overline B_R^X,
    \qquad
    \Phi(u_0,\mathcal S_r(u_0))=\mathcal S_r(u_0).
\end{equation*}
For an arbitrary initial datum, take $r=1+\norm{u_0}_{H^s}$ and then choose
$T$ as above. This proves existence on $[0,T]$, with $T$ depending only on
$\norm{u_0}_{H^s}$ (for fixed $s$ and $\delta$).

Next, we prove uniqueness in
$X_T^{s,\frac12+\delta}$. Let
$u,v\in X_T^{s,\frac12+\delta}$ be two solutions corresponding to the
same initial datum $u_0$. Fix
\begin{equation}
    0<T_0\leq T.
\end{equation}
For $0\leq t\leq T_0$, the Duhamel formula gives
\begin{equation}
\begin{aligned}
u(t)-v(t)
&=
-
\int_0^t
e^{-(t-t')\pt_x\Delta}
\pt_x
\prc{
    u^3(t')-v^3(t')
}
\,dt'.
\end{aligned}
\end{equation}
Since $\eta(t/T_0)=1$ for $|t|\leq T_0$, we may insert the cutoff
inside the nonlinear term and obtain, on $[0,T_0]$,
\begin{equation}
\begin{aligned}
u(t)-v(t)
&=
-
\eta\prc{\frac{t}{T_0}}
\int_0^t
e^{-(t-t')\pt_x\Delta}
\\
&\qquad {}\mathbin{\times}\pt_x
\left[
\begin{aligned}
    &\prc{\eta(t'/T_0)u(t')}^3\\
    &\quad-
    \prc{\eta(t'/T_0)v(t')}^3
\end{aligned}
\right]
\,dt'.
\end{aligned}
\end{equation}
For every $\varepsilon>0$, choose
$U,V,W\in X^{s,\frac12+\delta}$ such that
\begin{equation}
    U|_{(0,T_0)\times\rt}=u,
    \qquad
    V|_{(0,T_0)\times\rt}=v,
    \qquad
    W|_{(0,T_0)\times\rt}=u-v,
\end{equation}
and
\begin{equation}
\begin{aligned}
    \xnorm{U}{s}{\frac12+\delta}
    &\leq
    \norm{u}_{X_{T_0}^{s,\frac12+\delta}}
    +\varepsilon,
    \\
    \xnorm{V}{s}{\frac12+\delta}
    &\leq
    \norm{v}_{X_{T_0}^{s,\frac12+\delta}}
    +\varepsilon,
    \\
    \xnorm{W}{s}{\frac12+\delta}
    &\leq
    \norm{u-v}_{X_{T_0}^{s,\frac12+\delta}}
    +\varepsilon.
\end{aligned}
\end{equation}
Since, on $(0,T_0)\times\rt$,
\begin{equation}
    u^3-v^3
    =
    (u-v)(u^2+uv+v^2),
\end{equation}
the function
\begin{equation}
    W(U^2+UV+V^2)
\end{equation}
is an extension of $u^3-v^3$.

Using
\lemref{lem:duhamel-linear-estimate}
and
\lemref{lem:trilinear-bourgain},
we obtain
\begin{equation}
    \begin{aligned}
\norm{u-v}_{X_{T_0}^{s,\frac12+\delta}}
&\lesssim
T_0^{2\delta}
\Big(
    \xnorm{U}{s}{\frac12+\delta}^2
    +
    \xnorm{U}{s}{\frac12+\delta}
    \xnorm{V}{s}{\frac12+\delta}\\
    &\quad+
    \xnorm{V}{s}{\frac12+\delta}^2
\Big)
\xnorm{W}{s}{\frac12+\delta}.
\end{aligned}
\end{equation}
Letting $\varepsilon\downarrow0$, we obtain
\begin{equation}
\label{eq:uniqueness-local-contraction}
\begin{aligned}
\norm{u-v}_{X_{T_0}^{s,\frac12+\delta}}
&\lesssim
T_0^{2\delta}
\Big(
    \norm{u}_{X_{T_0}^{s,\frac12+\delta}}^2
    +
    \norm{u}_{X_{T_0}^{s,\frac12+\delta}}
    \norm{v}_{X_{T_0}^{s,\frac12+\delta}}
\\
&\qquad\qquad
    +
    \norm{v}_{X_{T_0}^{s,\frac12+\delta}}^2
\Big)
\norm{u-v}_{X_{T_0}^{s,\frac12+\delta}}.
\end{aligned}
\end{equation}
Since restriction can only decrease the restricted Bourgain norm,
\begin{equation}
    \norm{u}_{X_{T_0}^{s,\frac12+\delta}}
    \leq
    \norm{u}_{X_T^{s,\frac12+\delta}},
    \qquad
    \norm{v}_{X_{T_0}^{s,\frac12+\delta}}
    \leq
    \norm{v}_{X_T^{s,\frac12+\delta}}.
\end{equation}
Choosing $T_0>0$ sufficiently small so that the coefficient on the
right-hand side of
\eqref{eq:uniqueness-local-contraction}
is strictly smaller than one yields
\begin{equation}
    u=v
    \qquad
    \text{on }[0,T_0].
\end{equation}
Since $u(T_0)=v(T_0)$ and the Bourgain norms and estimates are invariant under
time translation, the same argument applies to the next interval. Repeating
this procedure finitely many
times proves uniqueness on $[0,T]$.

\smallskip
To prove that the data-to-solution map is real-analytic, define
\begin{equation}
    \mathcal F(u_0,u)
    :=
    u-\Phi(u_0,u).
\end{equation}
Thus
\begin{equation}
    \mathcal F:
    H^s(\rt)
    \times
    X_T^{s,\frac12+\delta}
    \longrightarrow
    X_T^{s,\frac12+\delta}.
\end{equation}
The estimates above show that $\Phi$ is continuous, affine in $u_0$, and a
cubic polynomial in $u$; hence $\mathcal F$ is real-analytic. Fix
$u_0\in B_r^H$ and let
\begin{equation}
    u=\mathcal S_r(u_0)
\end{equation}
be the solution obtained above. Then
\begin{equation}
    \mathcal F(u_0,u)=0.
\end{equation}
The Fr\'echet derivative of $\mathcal F$ with respect to the second
variable is
\begin{equation}
\begin{aligned}
D_u\mathcal F(u_0,u)[v]
&=
v-D_u\Phi(u_0,u)[v],
\end{aligned}
\end{equation}
where
\begin{equation}
\label{derivativeOfGamma}
\begin{aligned}
D_u\Phi(u_0,u)[v]
&=
-3
\eta\prc{\frac{t}{T}}
\int_0^t
e^{-(t-t')\pt_x\Delta}
\pt_x
\prc{
    u^2(t')v(t')
}
\,dt'.
\end{aligned}
\end{equation}
By \lemref{lem:duhamel-linear-estimate} and \lemref{lem:trilinear-bourgain},
we obtain
\begin{align}
\norm{
    D_u\Phi(u_0,u)[v]
}_{X_T^{s,\frac12+\delta}}
&\leq
3cT^{2\delta}
\norm{u}_{X_T^{s,\frac12+\delta}}^2
\norm{v}_{X_T^{s,\frac12+\delta}}
\nonumber\\
&\leq
3cT^{2\delta}R^2
\norm{v}_{X_T^{s,\frac12+\delta}}.
\end{align}
By the choice of $T$ in the contraction argument,
\begin{equation}
    3cT^{2\delta}R^2
    <
    \frac12.
\end{equation}
Hence
\begin{equation}
    \norm{
        D_u\Phi(u_0,u)
    }_{\operatorname{op}}
    <
    1.
\end{equation}
Therefore,
\begin{equation}
    D_u\mathcal F(u_0,u)
    =
    I-D_u\Phi(u_0,u)
\end{equation}
is invertible by
\thmref{Invertibility}.
Applying
\thmref{ImplicitFunctionTheorem}
to the map $\mathcal F$ at
\begin{equation}
    \prc{
        u_0,
        \mathcal S_r(u_0)
    },
\end{equation}
we obtain $\rho>0$ and a unique real-analytic map
\begin{equation}
    \widetilde{\mathcal S}:
    B_\rho(u_0)
    \longrightarrow
    X_T^{s,\frac12+\delta}
\end{equation}
such that
\begin{equation}
    \mathcal F
    \prc{
        v_0,
        \widetilde{\mathcal S}(v_0)
    }
    =
    0
\end{equation}
for every $v_0\in B_\rho(u_0)$. By the uniqueness already proved,
$\widetilde{\mathcal S}$ agrees with the data-to-solution map wherever both
are defined. Thus, the data-to-solution map is locally real-analytic at every
$u_0\in B_r^H$. Since
$X_{T'}^{s,\frac12+\delta}
\hookrightarrow C([0,T'];H^s(\rt))$
continuously, restriction to $[0,T']$ shows that the same
data-to-solution map is real-analytic into
\begin{equation*}
    C([0,T'];H^s(\rt))
    \cap
    X_{T'}^{s,\frac12+\delta}.
\end{equation*}
\end{proof}

\begin{remark}
It is also possible to prove analyticity by justifying the power series
expansion of the data-to-solution map; see Theorem~3 in
\cite{bejenaru2006sharp}.
\end{remark}

%% file: appendix.tex
\appendix
\section{Auxiliary results}
\subsection{Linear estimates in Bourgain's spaces}
In this section, we recall some standard results regarding Bourgain spaces; see \cite{ginibre1995probleme}. Let $\eta$ be the cutoff fixed in the proof of Theorem~\ref{MainResult}. The estimates below hold for any dispersion function $\omega$ and any spatial domain. The following lemma shows that the $X^{s,b}$ norm is stable under time localization.

\begin{lemma}
\label{Lem:XsbStability}
Let $s,b\in\R$ and $u\in X^{s,b}$. Then
$$\xnorm{\eta(t) u}{s}{b} \lesssim_{b} \xnorm{u}{s}{b}.$$
\end{lemma}

\begin{proof}
    Using Fourier inversion formula on $\eta$ and Minkowski's inequality, 

    \begin{align*}
        \xnorm{\eta(t) u}{s}{b} &= \xnorm{\prc{\frac{1}{2\pi}\int_\R \Fh{\eta} (\tau_0) e^{i t \tau_0 } d\tau_0 }u}{s}{b}\\
        & \le \int_\R \xnorm{\Fh{\eta}(\tau_0)e^{i t \tau_0} u}{s}{b} d\tau_0 \\
        &= \int_\R \abs{\Fh{\eta}(\tau_0)} \xnorm{e^{i t \tau_0} u}{s}{b} d \tau_0\\
        & = \int_\R \abs{\Fh{\eta}(\tau_0)} \norm{\inp{|(\xi,q)|}^s \inp{\tau - \omega(\xi,q)}^b \Fh{u}(\tau - \tau_0,\xi,q)}_{L^2_{\xi q} L^2_\tau} d\tau_0
    \end{align*}
Applying the change of variables $\tau' = \tau -\tau_0$ and $\inp{A+B}^b \lesssim_b \inp{A}^{|b|}\inp{B}^b$ gives 

\begin{align*}
    \xnorm{\eta(t) u}{s}{b} & \leq \int_\R \abs{\Fh{\eta}(\tau_0)} \norm{\inp{|(\xi,q)|}^s \inp{\tau'+\tau_0 - \omega(\xi,q)}^b \Fh{u}(\tau',\xi,q)}_{L^2_{\xi q} L^2_{\tau'}} d\tau_0\\
    & \le \int_\R \abs{\Fh{\eta}(\tau_0)} \inp{\tau_0}^{|b|} \norm{\inp{|(\xi,q)|}^s \inp{\tau' - \omega(\xi,q)}^b \Fh{u}(\tau',\xi,q)}_{L^2_{\xi q} L_{\tau'}^2} d\tau_0\\
    & \lesssim_b \xnorm{u}{s}{b}.
\end{align*}

\end{proof}

The following is Lemma~3.1 in \cite{ginibre1995probleme} with $T_0=1$ fixed. This choice is advantageous because, when using the contraction mapping principle in $X^{s,b}$ with $b>1/2$ and $T$ small, taking $T_0=1$ rather than $T_0=T$ allows us to work on a ball of fixed size independent of $T$.

\begin{lemma}[Homogeneous linear estimate] 
\label{Lem:homogeneous-linear-estimate}

Let $s,b\in\R$ and $u_0\in H^s(\R\times\T)$.
\begin{equation}
\label{lem:homogeneous-linear-estimate}
\xnorm{{\eta(t) e^{-t\pt_x\Delta}u_0}}{s}{b} \lesssim_b \Hnxy{u_0}{s}{xy}.
\end{equation}
\end{lemma}

\begin{proof}
    We start by calculating the space-time Fourier transform of $\eta\, \sg{-t} u_0$,

    \begin{align*}
        \F {\left({\eta(t) \sg{-t} u_0}\right)}(\tau, \xi,q) &= \frac{1}{(2\pi)^3} \int e^{- i (t \tau + (x,y) \cdot (\xi,q))} \eta(t) e^{i t \omega(\xi,q)} u_0(x,y) dx dy dt\\
        & = \brs{ \frac{1}{2\pi} \int_\R e^{- i t (\tau - \omega(\xi,q))} \eta(t) dt}\\
        &\quad \times \brs{\frac{1}{(2 \pi)^2} \int_{\R \times \T} e^{- i (x,y)\cdot (\xi,q)} u_0(x,y) dx dy }\\
        &= \Fh{\eta} (\tau - \omega(\xi,q)) \Fh{u_0}(\xi,q).
    \end{align*}
Now we can calculate the norm,

\begin{align*}
    \xnorm{\eta(t) \sg{-t} {u_0}}{s}{b} &= \norm{\inp{|(\xi,q)|}^s \inp{\tau - \omega(\xi,q)}^b \Fh{\eta}(\tau -\omega(\xi,q)) \Fh{u_0}(\xi,q)}_{L_{\xi q}^2 L^2_\tau}\\
    &= \Lptxy{  \Lptxy{\inp{\tau -\omega(\xi,q)}^b \Fh{\eta}(\tau - \omega(\xi,q))}{2}{\tau} \inp{|(\xi,q)|}^s \Fh{u_0}(\xi,q)  }{2}{\xi q}\\
    &= \norm{\eta}_{H^b_t} \norm{\inp{|(\xi,q)|}^s \Fh{u_0}}_{L^2_{\xi q}}\\
    &= \norm{\eta}_{H^b_t} \norm{u_0}_{H^{s}_{xy}}\\
    & \lesssim_b \norm{u_0}_{H^{s}_{xy}}.
\end{align*}

\end{proof}

\begin{lemma}[Inhomogeneous linear estimates]
\label{lem:duhamel-linear-estimate}
Let $s\in\R$,
$-\frac12<b'\leq0\leq b\leq b'+1$,
$0<T\leq1$, and $F\in X^{s,b'}$. Then,

\begin{equation}
    \label{eq:duhamel-linear-estimate}
    \norm{\eta(t/T) \int_0^t e^{-(t-t')\pt_x \Delta} F(t') dt'}_{X^{s,b}} \le c\, T^{1-b+b'} \norm{F}_{X^{s,b'}},
\end{equation}

where $c=c(\eta,b,b')>0$ is independent of $T$ and $F$.
\end{lemma}

\begin{proof} 
See Lemma 2.1 (ii) in \cite{ginibre1997cauchy} and Lemma 3.2 in \cite{ginibre1995probleme}.
\end{proof}

\begin{lemma}[Time Localization Estimate]
\label{Lem:XbsLocalImbedding}

Let $0<T<1$, $s\in\R$, $-\frac12<b'\leq b<\frac12$, and
$u\in X^{s,b}$. Then

\begin{equation}
    \norm{\eta(t/T)u}_{X^{s,b'}} \lesssim_{b,b'} T^{b-b'} \norm{u}_{X^{s,b}}.
\end{equation}
    
\end{lemma}
\begin{proof}
Set
\begin{equation}
    U(t)
    :=
    e^{t\pt_x\Delta}u(t).
\end{equation}
By the definition of the Bourgain norm,
\begin{equation}
    \xnorm{u}{s}{b}
    =
    \norm{U}_{H_t^bH_{xy}^s},
\end{equation}
and, since multiplication by $\eta(t/T)$ commutes with the spatial
group,
\begin{equation}
    \xnorm{\eta(t/T)u}{s}{b'}
    =
    \norm{
        \eta(t/T)U
    }_{H_t^{b'}H_{xy}^s}.
\end{equation}
The standard one-dimensional time-localization estimate
\begin{equation}
    \norm{
        \eta(t/T)F
    }_{H_t^{b'}(\R;H)}
    \lesssim_{b,b'}
    T^{b-b'}
    \norm{F}_{H_t^b(\R;H)},
\end{equation}
valid for every Hilbert space $H$ and
\begin{equation}
    -\frac12<b'\leq b<\frac12,
    \qquad
    0<T<1,
\end{equation}
therefore gives, with $H=H^s(\rt)$,
\begin{equation}
    \xnorm{\eta(t/T)u}{s}{b'}
    \lesssim_{b,b'}
    T^{b-b'}
    \xnorm{u}{s}{b}.
\end{equation}
The one-dimensional estimate is standard; see
Lemma~2.11 in \cite{tao2006nonlinear} and
Lemma~3.11 in \cite{erdougan2016dispersive}.
\end{proof}

\subsection{Analytic maps and the implicit function theorem}
\label{sec:appendix-analytic-maps}
In this section, we present the mathematical framework for analytic maps and the implicit function theorem in Banach spaces. The following definitions and theorems are standard; see \cite{herr2006well} for details.

\begin{definition}
Let $X,Y$ be Banach spaces over $\R$ or $\C$, and let $U\subset X$
be open. We denote by $C^1(U,Y)$ the set of all Fr\'echet
differentiable maps $F:U\to Y$ such that
\begin{equation}
    F':U\to L(X,Y)
\end{equation}
is continuous. Inductively, we define $C^k(U,Y)$ for
$k=2,\ldots,\infty$.
\end{definition}
\begin{definition}
Let $X,Y$ be Banach spaces over $\R$, let $U\subset X$ be open, and
let $F:U\to Y$. We say that $F$ is real-analytic if, for every
$u\in U$, there exists $r>0$ such that
\begin{equation}
    B_r(u)\subset U,
\end{equation}
and there exist continuous $k$-linear maps
\begin{equation}
    L_k:X^k\to Y,
    \qquad
    k\in\N_0,
\end{equation}
such that, writing
\begin{equation}
    L^{(k)}(x):=L_k(x,\ldots,x),
\end{equation}
one has
\begin{equation}
    F(x)
    =
    \sum_{k=0}^{\infty}
    L^{(k)}(x-u),
    \qquad
    x\in B_r(u),
\end{equation}
with convergence in $Y$.
\end{definition}

\begin{theorem}[Proposition 1.4.1 in \cite{herr2006well}]
Let $X,Y$ be Banach spaces over $\R$ or $\C$, let $U\subset X$ be open,
and let $F:U\to Y$ be differentiable. Assume that for every $u\in U$ there
exists $\epsilon>0$ such that $B_\epsilon(u)\subset U$ and
$F':B_\epsilon(u)\to L(X,Y)$ is bounded. Then
$F\big|_{B_\epsilon(u)}$ is Lipschitz continuous. In particular, if
$F\in C^1(U,Y)$, then $F$ is locally Lipschitz continuous.
\end{theorem}
\begin{theorem}[Proposition 1.4.5 in \cite{herr2006well}]
Let $X,Y,Z$ be Banach spaces over $\R$, let $U\subset X$ and
$V\subset Y$ be open, and let
\begin{equation}
    F:U\to V,
    \qquad
    G:V\to Z
\end{equation}
be real-analytic maps. Then
\begin{equation}
    G\circ F:U\to Z
\end{equation}
is real-analytic.
\end{theorem}

\begin{theorem}[Neumann-series criterion]
\label{Invertibility}
Let $X$ be a Banach space and let $T\in L(X,X)$ satisfy
\begin{equation}
    \norm{T}_{L(X,X)}<1.
\end{equation}
Then $I-T$ is boundedly invertible on $X$, and
\begin{equation}
    (I-T)^{-1}
    =
    \sum_{n=0}^{\infty}T^n
\end{equation}
with convergence in $L(X,X)$.
\end{theorem}
\begin{theorem}[\textbf{Implicit Function Theorem}]
\label{ImplicitFunctionTheorem}
Let $k\in\N$, let $X,Y,Z$ be Banach spaces over $\R$, and let
$U\subset X$ and $V\subset Y$ be open neighborhoods of $x_0$ and
$y_0$, respectively. Let
\begin{equation}
    F:U\times V\to Z
\end{equation}
be of class $C^k$ (respectively, real-analytic), and assume that
\begin{equation}
    F(x_0,y_0)=0.
\end{equation}
Assume moreover that
\begin{equation}
    D_yF(x_0,y_0):Y\to Z
\end{equation}
is a bounded linear isomorphism, with
\begin{equation}
    \big(D_yF(x_0,y_0)\big)^{-1}\in L(Z,Y).
\end{equation}
Then there exist neighborhoods $B_r(x_0)\subset U$ and
$B_R(y_0)\subset V$ and a unique map of class $C^k$ (respectively,
real-analytic)
\begin{equation}
    G:B_r(x_0)\to B_R(y_0)
\end{equation}
such that
\begin{equation}
    G(x_0)=y_0
\end{equation}
and
\begin{equation}
    F(x,G(x))=0
\end{equation}
for every $x\in B_r(x_0)$.
\end{theorem}

%% file: acknowledgment.tex
\section*{Acknowledgments}
\noindent
The author would like to acknowledge Professors Kenji Nakanishi and Nobu Kishimoto for their significant guidance throughout much of this work, which was conducted at the Research Institute for Mathematical Sciences (RIMS), Kyoto University, from 2021 to 2023. The author gratefully acknowledges support from NSF grant DMS-2009800 through Professor Thomas Chen.